\documentclass[a4paper,11pt]{amsart}
\usepackage{amssymb}
\usepackage{mathrsfs}
\usepackage{amsmath,amssymb,amsthm,latexsym,amscd,mathrsfs}
\usepackage{indentfirst}
\usepackage{stmaryrd}
\usepackage{graphicx}
 \newcommand{\ROM}[1]{\mathrm{\uppercase\expandafter{\romannumeral#1}}}
  \theoremstyle{definition}

 \newtheorem{thm}{Theorem}[section]
 \newtheorem{lem}{Lemma}[section]

 \newtheorem{rem}{Remark}[section]
 \newtheorem{prop}{Proposition}[section]

\newtheorem{ack}{Acknowledgements}   
  \makeatletter

  \numberwithin{equation}{section}
  \allowdisplaybreaks
\makeatother
\title[Isoparametric foliation and a problem of Besse]{\textbf{Isoparametric foliation and a problem of Besse on generalizations of Einstein condition}}
\author[Z.Z.Tang]{Zizhou Tang}\address{School of Mathematical Sciences, Laboratory of Mathematics and Complex Systems, Beijing Normal
University, Beijing 100875, China}\email{zztang@bnu.edu.cn}
\author[W. J. Yan]{Wenjiao Yan$^{\dag}$}
\address{School of Mathematical Sciences, Laboratory of Mathematics and Complex Systems, Beijing Normal
University, Beijing 100875, China} \email{wjyan@bnu.edu.cn}
\thanks {$^{\dag}$ the corresponding author}

\date{}

\begin{document}

\maketitle

\begin{abstract}
 The focal sets of isoparametric
hypersurfaces in spheres with $g=4$ are all Willmore submanifolds,
being minimal but mostly non-Einstein (\cite{TY1}, \cite{QTY}).
Inspired by A.Gray's view, the present paper shows that, these focal
sets are all $\mathcal{A}$-manifolds but rarely Ricci parallel,
except possibly for the only unclassified case. As a byproduct, it
gives infinitely many simply-connected examples to the problem 16.56
(i) of Besse concerning generalizations of the Einstein condition.
\end{abstract}

\section{\textbf{Introduction}}
The Riemannian manifolds with constant Ricci curvatures (the
Einstein condition) and those with constant scalar curvatures are
two important classes of Riemannian manifolds. We denote them by
$\mathcal{E}$ and $\mathcal{S}$, respectively. Then there comes
apparently a class of manifolds with parallel Ricci tensor, denoted by
$\mathcal{P}$, lies between $\mathcal{E}$ and $\mathcal{S}$. As
further generalizations of the Einstein condition, A. Gray
(\cite{Gra}) introduced two significant classes $\mathcal{A}$ and
$\mathcal{B}$ defined as follows, in which the Ricci tensor $\rho$
is cyclic parallel and a Codazzi tensor, respectively:
\begin{eqnarray}\label{A manifold}
  \mathcal{A} &:& \nabla_i\rho_{jk}+\nabla_j\rho_{ki}+\nabla_k\rho_{ij}=0 \\
  \mathcal{B} &:& \nabla_i\rho_{jk}-\nabla_j\rho_{ik}=0.\nonumber
\end{eqnarray}
These two classes $\mathcal{A}$ and $\mathcal{B}$ are investigated
extensively since then. In  view of the  second Bianchi identity,
the class $\mathcal{B}$ coincides with those having harmonic
curvatures. Gray also showed that the following inclusions exist
between the various classes:
\begin{equation*}
  \mathcal{E}\subset\mathcal{P}=\mathcal{A\cap B}
  _{{\rotatebox{-45}{$\subset \rotatebox{45}{$\mathcal{B}$}$}}}
    ^{{\rotatebox{45}{$ \subset {\rotatebox{-45}{$\mathcal{A}$}}$}}}
    {_{^{\rotatebox{45}{$\subset~~~~$}}}^{\rotatebox{-45}{$\subset$}}\mathcal{A}\cup\mathcal{B}}\subset \mathcal{S}
\end{equation*}
and proved that $\mathcal{A}$ and $\mathcal{B}$ are the only classes between $\mathcal{P}$ and $\mathcal{S}$ from the
view of group representations.

\cite{TY1} and \cite{QTY} provide many new examples of Willmore
submanifolds in spheres via isoparametric foliation. More precisely,
the focal submanifolds of isoparametric hypersurfaces in spheres
with four distinct principal curvatures (the most complicated and
abundant case) are all Willmore submanifolds in spheres. Since the
focal submanifolds are minimal in spheres, in contrast with the
well-known fact that the Einstein manifolds minimally immersed in
spheres are Willmore, \cite{QTY} also determined which of these
focal submanifolds are Einstein. A further and natural question
arises: \emph{are they Ricci parallel, $\mathcal{A}$-manifolds, or
$\mathcal{B}$-manifolds ? }The present paper aims at an answer to
this question.

To state clearly the results, we first need a short review of the
isoparametric foliation.

Researches on classifications and applications of isoparametric
foliation in spheres have been quite active recently (for classifications,
see \cite{CCJ}, \cite{Miy}, \cite{Chi}; for applications, see for
example, \cite{GR}, \cite{GX}, \cite{QTY}, \cite{TXY}, \cite{TY1}, \cite{TY2}). As is well
known, an isoparametric hypersurface $M^n$ in $S^{n+1}(1)$ is a
hypersurface with constant principal curvatures. Let $g$ be the
number of distinct principal curvatures with multiplicity $m_i$
($i=1,...,g$). According to \cite{Mun}, $g$ can be only $1$, $2$,
$3$, $4$ or $6$, and $m_i=m_{i+2}$ (subscripts mod $g$). When $g=1,
2, 3$ and $6$, the classification for isoparametric hypersurfaces
are accomplished; when $g=4$, all isoparametric hypersurfaces are of
OT-FKM type (defined later), or of homogeneous type with $(m_1,
m_2)=(2, 2), (4, 5)$ except possibly for the case with $(m_1,
m_2)=(7, 8)$ (cf. \cite{Chi}).

In fact, an isoparametric hypersurfaces in $S^{n+1}(1)$ always comes
as a family of parallel hypersurfaces, which are level hypersurfaces
(isoparametric foliation) of an isoparametric function $f$, that is,
a function on $S^{n+1}(1)$ satisfying
\begin{equation}\label{ab}
\left\{ \begin{array}{ll}
|\nabla f|^2= b(f),\\
~~~~\Delta f~~=a(f),
\end{array}\right.
\end{equation}
where $\nabla f$ and $\Delta f$ are the gradient and Laplacian of
$f$ on $S^{n+1}(1)$, $b$ and $a$ smooth and continuous functions on
$\mathbb{R}$, respectively. The two singular sets of $f$ are called
the focal sets (submanifolds), denoted respectively by $M_{1}$ and
$M_{2}$, being actually minimal submanifolds of $S^{n+1}(1)$ with
codimensions $m_1+1$ and $m_2+1$ (cf. \cite{CR}).

Now we recall the construction of isoparametric functions of OT-FKM
type, constructed by Ferus, Karcher and M\"{u}nzner (\cite{FKM}),
following \cite{OT}. For a symmetric Clifford system
$\{P_0,\cdots,P_m\}$ on $\mathbb{R}^{2l}$, \emph{i.e.}
$P_{\alpha}$'s are symmetric matrices satisfying
$P_{\alpha}P_{\beta}+P_{\beta}P_{\alpha}=2\delta_{\alpha\beta}I_{2l}$,
a homogeneous polynomial $F$ of degree $4$ on $\mathbb{R}^{2l}$ is
defined as:
\begin{eqnarray}\label{FKM isop. poly.}
&&F(x) = |x|^4 - 2\displaystyle\sum_{\alpha = 0}^{m}{\langle
P_{\alpha}x,x\rangle^2}.
\end{eqnarray}
It is easy to verify that $f=F|_{S^{2l-1}}$ is an isoparametric function
on $S^{2l-1}$. The focal submanifolds $M_1=f^{-1}(1)$, $M_2=f^{-1}(-1)$, and the
multiplicity pair is $(m_1, m_2)=(m, l-m-1)$, provided
$m>0$ and $l-m-1>0$, where $l = k\delta(m)$ $(k=1,2,3,\cdots)$, $\delta(m)$ is the dimension
of an irreducible module of the Clifford algebra $\mathcal{C}_{m-1}$.

It was observed by \cite{KN} that the isoparametric hypersurfaces
are $\mathcal{A}$-manifolds only when $g\leq 3$, and Ricci parallel
only when $g\leq 2$. The present paper will study in-depth the focal
submanifolds. As one of the main results, we prove

\begin{thm}\label{class A}
All the focal submanifolds of isoparametric hypersurfaces in spheres
with $g=4$ are $\mathcal{A}$-manifolds, except possibly for the only
unclassified case with $(m_1, m_2)=(7, 8)$.
\end{thm}

From now on, we shall concentrate on the focal submanifolds $M_i$
($i=1, 2$) not in the unclassified case $(m_1, m_2)=(7, 8)$. From
Theorem \ref{class A} and the relation
 $\mathcal{P}$=$\mathcal{A\cap B}$
 it follows that
\begin{center}
   $M_i$ $\in\mathcal{P}$ $\Longleftrightarrow$ $M_i\in\mathcal{B}$, ($i=1, 2$).
\end{center}
Thus we are left to investigate which $M_i$ is Ricci parallel.

The following theorem achieves a complete answer to this question.

\begin{thm}\label{Ricci parallel}
For the focal submanifolds of isoparametric hypersurfaces in spheres
with $g=4$, we have
\begin{itemize}
\item [(i)] The $M_1$ of OT-FKM type is Ricci parallel if and only if
$(m_1,m_2)=(2,1), (6,1)$, or it is diffeomorphic to $Sp(2)$ in the
homogeneous case with $(m_1,m_2)=(4,3)$; while the $M_2$ of OT-FKM
type is Ricci parallel if and only if $(m_1,m_2)=(1,k)$.

\item[(ii)]For $(m_1,m_2)=(2,2)$, the one
           diffeomorphic to $\widetilde{G}_2(\mathbb{R}^5)$ is Ricci parallel, while the other
           diffeomorphic to $\mathbb{C}P^3$ is not.
\item[(iii)]For $(m_1,m_2)=(4,5)$, both are not Ricci parallel.

\end{itemize}
\end{thm}

\begin{rem}
As we mentioned in \cite{QTY}, the only Einstein ones among the known focal
submanifolds are actually the focal submanifold $M_1$ of OT-FKM type diffeomorphic
to $Sp(2)$ in the homogeneous case with $(m_1,m_2)=(4,3)$ and the focal submanifold
diffeomorphic to $\widetilde{G}_2(\mathbb{R}^5)$ with $(m_1,m_2)=(2,2)$.
\end{rem}

It is well known that the D'Atri spaces (Riemannian manifolds with
volume preserving geodesic symmetries) belong to the class
$\mathcal{A}$ (cf. \cite{Bes}, pp.450). So the examples of
$\mathcal{A}$-manifolds are not rare in the literature, but mostly
are (locally) homogeneous.

In this regard, Besse (\cite{Bes}, 16.56(i), pp.451) posed the
following problem as one of ``some open problems" : \emph{Find
examples of $\mathcal{A}$-manifolds, which are neither locally
homogeneous, nor locally isometric to Riemannian products and have
non-parallel Ricci tensor}.

To find examples for this problem, \cite{Jel} and \cite{PT}
constructed $\mathcal{A}$-manifolds on $S^1$-bundles over locally
non-homogeneous K\"{a}hler-Einstein manifolds, and on $S^1$-bundles
over a $K3$ surface, from defining Riemannian submersion metric on
the $S^1$-bundles. But in some sense, their examples are not so satisfying,
as they are not simply-connected, and the metrics are not
natural enough.

On the ground of Theorem \ref{class A} and \ref{Ricci parallel}, in
virtue of the following two propositions, we find a series of
natural, simply-connected examples for this open problem of Besse.

\begin{prop}\label{Rie product}
The focal submanifolds of isoparametric hypersurfaces in spheres with $g=4$
and $m_1, m_2>1$ are not Riemannian products.
\end{prop}

\begin{prop}\label{locally homogeneous}
The focal submanifolds $M_1$ of OT-FKM type with $(m_1, m_2)=(3, 4k)$
are not intrinsically homogeneous.
\end{prop}

\begin{rem}
By Morse theory, one sees that if $m_1>1$ (resp. $m_2>1$), the focal
submanifold $M_2$ (resp. $M_1$) is simply-connected (cf.
\cite{Tan}). Combining the two propositions above with Theorem
\ref{class A} and \ref{Ricci parallel}, we conclude that the focal
submanifolds $M_1$ of OT-FKM type with $(m_1, m_2)=(3, 4k)$ are
simply-connected $\mathcal{A}$-manifolds with non-parallel Ricci
tensor, which are minimal submanifolds in spheres, but neither
locally homogeneous, nor locally isometric to Riemannian products.
Much more examples to the problem of Besse can be obtained in this
way, however we shall not go into the details in this paper.
\end{rem}


\section{\textbf{$\mathcal{A}$-manifolds}}

We begin this section with displaying an equivalent condition of the
definition (\ref{A manifold}) for $\mathcal{A}$-manifold, that is
\begin{equation}\label{equi A}
  (\nabla_X\rho)(X,X)=0, \mathrm{~for ~ any ~tangent~ vector ~}X.
\end{equation}
Based on the known classification results of the isoparametric
hypersurfaces in spheres with four distinct principal curvatures, we
will divide the proof of Theorem \ref{class A} into three parts: the
OT-FKM type, the homogeneous cases with $(m_1, m_2)=(2, 2)$ and $(4,
5)$.

\subsection{OT-FKM type}

\subsubsection{\textbf{$M_1$ of OT-FKM type}.}
According to the definition (\ref{FKM isop. poly.}), the focal submanifold $M_1$ of OT-FKM type
can be written as:
$$M_1=\{x\in S^{2l-1}~|~\langle P_0x, x\rangle=\langle P_1x, x\rangle=\cdots=\langle P_mx, x\rangle=0\}.$$
Note that $\dim M_1=2l-m-2$. As pointed out by \cite{FKM}, the
normal space in $S^{2l-1}$ at $x\in M_1$ is
\begin{equation*}
T^{\perp}_xM_1=\{~Px~|~P\in\mathbb{R}\Sigma (P_0,...,P_m)~\},
\end{equation*}
where $\Sigma (P_0,...,P_m)$ is the unit sphere in $Span\{P_0,...,P_m\}$, which is called \emph{the Clifford sphere}
determined by the system $\{P_0,...,P_m\}$.

For the normal vector $\xi_{\alpha}=P_{\alpha}x$, $\alpha=0,...,m$,
denote $A_{\alpha}=:A_{\xi_{\alpha}}$ the shape operator
corresponding to $\xi_{\alpha}$. Then for any $X,Y \in T_xM_1$, the
Ricci tensor $\rho(X, Y)$ is given by (cf. \cite{TY1})
\begin{equation*}
\rho(X, Y)=2(l-m-2)\langle X, Y\rangle + \sum_{\alpha, \beta =0, \alpha \neq \beta}^{m}\langle X, P_{\alpha}P_{\beta}x\rangle\langle Y, P_{\alpha}P_{\beta}x\rangle.
\end{equation*}
As the metric tensor is parallel, we need only to focus on the tensor
$$\sigma(X, Y)=:\sum_{\alpha, \beta =0, \alpha \neq \beta}^{m}\langle X, P_{\alpha}P_{\beta}x\rangle\langle Y, P_{\alpha}P_{\beta}x\rangle.$$
A direct calculation implies
\begin{eqnarray}\label{M1 class A}
&&(\nabla_Z\sigma)(X, Y) \nonumber\\
&=& Z\sigma(X, Y)-\sigma(\nabla_ZX, Y)- \sigma(X, \nabla_ZY)\nonumber\\
   &=&  \sum_{\alpha, \beta =0, \alpha \neq \beta}^{m}\langle X, ~\nabla_Z(P_{\alpha}P_{\beta}x)\rangle\langle Y, ~P_{\alpha}P_{\beta}x~\rangle+\langle X, ~P_{\alpha}P_{\beta}x~\rangle\langle Y,~ \nabla_Z(P_{\alpha}P_{\beta}x)\rangle\nonumber\\
   &=& \sum_{\alpha, \beta =0, \alpha \neq \beta}^{m} \langle X,~ D_Z(P_{\alpha}P_{\beta}x)\rangle\langle Y, ~P_{\alpha}P_{\beta}x\rangle+\langle X, ~ P_{\alpha}P_{\beta}x\rangle\langle Y,~ D_Z(P_{\alpha}P_{\beta}x)\rangle\\
   &=& -\sum_{\alpha, \beta =0, \alpha \neq \beta}^{m}\langle Z, ~~\langle Y, P_{\alpha}P_{\beta}x\rangle P_{\alpha}P_{\beta}X+\langle X, P_{\alpha}P_{\beta}x\rangle P_{\alpha}P_{\beta}Y\rangle,\nonumber
\end{eqnarray}
where $\nabla$ and $D$ are the Levi-Civita connections on $M_1$ and
$\mathbb{R}^{2l}$, respectively.

Apparently, taking $X=Y=Z$, (\ref{M1 class A}) leads directly to
$(\nabla_X\sigma)(X, X)=0$, equivalently, $M_1$ of OT-FKM type is an
$\mathcal{A}$-manifold, as we desired.

\subsubsection{\textbf{$M_2$ of OT-FKM type}.}
Following \cite{FKM}, we see
that the focal submanifold
\begin{eqnarray*}
M_2 &=&F^{-1}(-1)\cap S^{2l-1}\\
&=& \{x\in S^{2l-1} |\  \mathrm{there \ exists} \ P \in \Sigma(P_0,\cdots,P_m) \ with\  Px=x\}.
\end{eqnarray*}
Observe that for any $P\in \Sigma(P_0,\cdots,P_m)$, its eigenvalues
must be $\pm 1$, with equal multiplicity. Thus $\mathbb{R}^{2l}$ can
be decomposed as a direct sum of the corresponding eigenspaces
$E_+(P)$ and $E_-(P)$.

Given $x \in M_2$ and $P\in \Sigma(P_0,\cdots,P_m)$ with $Px=x$. Define
$$\Sigma_P:=\{Q\in \Sigma(P_0,\cdots,P_m)|\  \langle P, Q \rangle:=\frac{1}{2l}\mathrm{Trace}(PQ)=0 \}, $$
which is the equatorial sphere of $\Sigma(P_0,\cdots,P_m)$ orthogonal to $P$.
In this way, there exists a decomposition of the tangent space $T_xM_2$
with respect to the eigenspaces of the shape operator.

\noindent \textbf{Lemma (\cite{FKM})}\,\, {\itshape The principal
curvatures of the shape operator $A_{\eta}$ with respect to any unit
normal vector $\eta\in T^{\perp}_xM_2$ are $0, 1$, and $-1$, with
the corresponding eigenspaces $Ker(A_{\eta})$, $E_+(A_{\eta})$,
$E_-(A_{\eta})$ as follows:
\begin{eqnarray}\label{eigenspacesM-}
Ker(A_{\eta})&=& \{v\in E_+(P)| \ v \bot x,\  v\bot \Sigma_P \eta\},\nonumber\\
E_+(A_{\eta})&=& \mathbb{R}\Sigma_P (x+\eta),\\
E_-(A_{\eta})&=& \mathbb{R}\Sigma_P (x-\eta).\nonumber
\end{eqnarray}Moreover, $$\dim Ker(A_{\eta})=l-m-1,\ \dim E_+(A_{\eta})=\dim E_-(A_{\eta})=m.$$
}

Let's now choose $\eta_1, \eta_2, \cdots, \eta_{l-m}$ as an
orthonormal basis of $T^{\perp}_xM_2$ in $S^{2l-1}$. Denote
$A_{\alpha}=:A_{\eta_{\alpha}}$. Then the minimality of $M_2$ in
$S^{2l-1}$ leads us to the following expression of the Ricci tensor
with respect to $X, Y\in T_xM_2$:
\begin{equation}\label{rho(X,Y) of M2}
 \rho(X, Y)=(l+m-2)\langle X, Y\rangle-\sum_{\alpha=1}^{l-m}\langle A_{\alpha}X, A_{\alpha}Y \rangle.
\end{equation}
Again, we just need to deal with the tensor $\displaystyle \tau(X,
Y)=:\sum_{\alpha=1}^{l-m}\langle A_{\alpha}X, A_{\alpha}Y \rangle$.

For this purpose, we make some preparation. In order to facilitate
the expression, we denote $Q_0=:P$. Then we can extend it to such a
symmetric Clifford system $\{Q_0, Q_1, \cdots , Q_m\}$ with
$Q_i\in\Sigma_P$ $(i\geq 1)$ that $\Sigma(Q_0, Q_1, \cdots,
Q_m)=\Sigma(P_0, P_1, \cdots, P_m)$. Using the previous lemma, it is
not difficult to find the following:

\begin{lem}
Given $i$, $1\leq i\leq m$, the unit vectors
\begin{equation}\label{o.n.b of TxM2}
\{
Q_i\eta_1,\cdots,Q_i\eta_{l-m},~~Q_1x,\cdots,Q_mx,~~Q_iQ_1x,\cdots,\widehat{Q_iQ_{i}x},\cdots,Q_iQ_mx\}.
\end{equation}
constitute an orthonormal basis of $T_xM_2$.
\end{lem}

Observe that by (\ref{eigenspacesM-}), we can decompose $A_{\alpha}X$ as
$$A_{\alpha}X=\sum_{i=1}^m\langle X, Q_ix\rangle Q_i\eta_{\alpha} + \langle X, Q_i\eta_{\alpha}\rangle Q_ix.$$
Then a direct verification by using (\ref{o.n.b of TxM2}) shows that
\begin{eqnarray}\label{AX AY}
  \tau(X, Y)
   &=& \sum_{\alpha=1}^{l-m}\sum_{i=1}^m\langle X, Q_ix\rangle\langle Y, Q_ix\rangle+\langle X, Q_i\eta_{\alpha}\rangle\langle Y, Q_i\eta_{\alpha}\rangle \nonumber\\
   &=& \sum_{\alpha=1}^{l-m}\sum_{k=0}^m \langle X, P_kx\rangle\langle Y, P_kx\rangle\\
   &&+ \sum_{i=1}^m\big\{~\langle X, Y\rangle-\sum_{j=1}^m\langle X, Q_jx\rangle\langle Y, Q_jx\rangle-\sum_{j=1,j\neq i}^m\langle X, Q_iQ_jx\rangle\langle Y, Q_iQ_jx\rangle\big\}\nonumber\\
   &=& m\langle X, Y\rangle+(l-2m)\sum_{k=0}^m \langle X, P_kx\rangle\langle Y, P_kx\rangle-\sum_{i,j=1,i\neq j}^m\langle X, Q_iQ_jx\rangle\langle Y, Q_iQ_jx\rangle\nonumber
\end{eqnarray}
Define
\begin{equation}\label{V W}
  V(X, Y)=:\sum_{k=0}^m \langle X, P_kx\rangle\langle Y, P_kx\rangle,\quad W(X, Y)=:\sum_{i,j=1,i\neq j}^m\langle X, Q_iQ_jx\rangle\langle Y, Q_iQ_jx\rangle.
\end{equation}
Then decomposing $X$ with respect to the orthonormal basis
(\ref{o.n.b of TxM2}), we see that the tensor $V$ is cyclic
parallel, since
\begin{eqnarray}\label{V class A}
  (\nabla_XV)(X, X) &=& 2\sum_{k=0}^m\langle X, P_kx\rangle\langle X, P_kX\rangle  \nonumber\\
   &=& 2\sum_{i=1}^m\langle X, Q_ix\rangle\langle X, Q_iX\rangle  \\
   &=& 4\sum_{i,j=1,i\neq j}^m \langle X, Q_ix\rangle\langle X, Q_jx\rangle\langle X, Q_iQ_jx\rangle\nonumber\\
   &=& 0. \nonumber
\end{eqnarray}
As for the tensor $W$, we can rewrite it as
\begin{equation}\label{W}
  W(X, Y) = \sum_{k,s=0, k\neq s}^m \langle X, P_kP_sx\rangle\langle Y, P_kP_sx\rangle-2V(X, Y)
\end{equation}
Thus it is easy to see that
\begin{equation}\label{W class A}
(\nabla_XW)(X, X) = 2\sum_{k,s=0, k\neq s}^m \langle X, P_kP_sx\rangle\langle X, P_kP_sX\rangle -2(\nabla_XV)(X, X)=0.
\end{equation}
At last, combining (\ref{rho(X,Y) of M2}), (\ref{AX AY}), (\ref{V class A}) with (\ref{W class A}), we
arrive at the conclusion that the focal submanifold $M_2$ of OT-FKM type is an $\mathcal{A}$-manifold, as desired.

\subsection{the homogeneous case}

It is well known that a homogeneous (isoparametric) hypersurface in
$S^{n+1}(1)$ can be characterized as a principal orbit of the
isotropy representation of some rank two symmetric space $G/K$,
while focal submanifolds correspond to the singular orbits (cf.
\cite{HL}). Denote by $\mathcal{G}$ and $\mathfrak{k}$ the Lie
algebras of $G$ and $K$, respectively. Then one has the following
Cartan decomposition
\begin{equation*}
  \mathcal{G}\cong \mathfrak{k}\oplus\mathfrak{p}.
\end{equation*}
Let $\langle\cdot ,\cdot\rangle$ be the usual $Ad(K)$-invariant
inner product on $\mathcal{G}$ that is induced from the Killing form
and the Cartan involution of $\mathcal{G}$. Following \cite{BCO},
let $z_0\in\mathfrak{p}$ be a unit vector and $M=Ad(K)\cdot z_0$ the
corresponding adjoint orbit included in the unit sphere. This leads
to a reductive decomposition of $\mathfrak{k}$ at $z_0$:
$$\mathfrak{k}=\mathfrak{k}_{z_0}\oplus\mathfrak{m},$$
where $\mathfrak{k}_{z_0}=\{ Y\in \mathfrak{k}~|~[Y, z_0]=0\}$ is
the isotropy subalgebra at $z_0$, and $\mathfrak{m}$ is the
orthogonal complement with respect to $\langle\cdot ,\cdot\rangle$
of $\mathfrak{k}_{z_0}$ in $\mathfrak{k}$. The tangent space and
normal space of $M$ in $\mathfrak{p}$ at $z_0$ are given by
\begin{equation*}
  T_{z_0}M=[\mathfrak{m}, z_0], \quad T_{z_0}^{\perp}M=\{~\xi\in\mathfrak{p}~|~[\xi, z_0]=0\},
\end{equation*}
while the shape operator with respect to $\xi$ is
\begin{equation*}
  A_{\xi}[m, z_0]=-[m, \xi]^{\top}, ~\mathrm{for} ~m\in
  \mathfrak{m},
\end{equation*}
where $(\cdot)^{\top}$ denotes the orthogonal projection to
$T_{z_0}M$. We prepare the following lemma whose proof is omitted.
\begin{lem}
Given $m, \tilde{m}\in \mathfrak{m}$, the Levi-Civita connection on
$M$ is stated as
$$\nabla_{[m, z]}[\tilde{m}, z]=[\tilde{m}, [m, z]]^{\top}, ~z\in M.$$
\end{lem}

Now let $M^n=Ad(K)\cdot z_0$ be a singular orbit, so that it is a
minimal submanifold in the unit sphere $S^{n+p}$. Choose
$\xi_1,\cdots,\xi_p$ as a unit normal basis. Similar as (2.4), to
verify the condition (2.1), we need only deal with the tensor
$\displaystyle \tau(X, Y)=\sum_{\alpha=1}^p\langle A_{\alpha}X,
A_{\alpha}Y\rangle$, for which we have
\begin{equation*}\label{u class A}
\frac{1}{2}(\nabla_X\tau)(X,X) = \sum_{\alpha=1}^p\langle~
\nabla_X(A_{\alpha}X)-A_{\alpha}(\nabla_XX), A_{\alpha}X\rangle.
\end{equation*}
Given a tangent vector at $z_0$, say $[m, z_0]$ for some $m\in
\mathfrak{m}$, we extend it to a tangent vector field $X$ by
$X(z)=[m, z]$ for $z\in M$. We have the following two equations and define $m'$
uniquely by the first one:
$$ \nabla_XX~|_{z_0}=[m, [m, z_0]]^{\top}=[m', z_0],\quad A_{\alpha}X=-[m, \xi_{\alpha}]^{\top}.$$
Let $\gamma(t)=exp(tm)\cdot z_0\in M$ be a curve so that
$\gamma(0)=z_0$, $\gamma'(0)=X(z_0)$. Clearly, for any
$\alpha=1,\cdots, p$, the unit normal vector $\xi_{\alpha}$ at $z_0$
can be extended along the curve $\gamma(t)$ to
$\xi_{\alpha}(t)=exp(tm)\cdot \xi_{\alpha}$. Then it is easy to see
that $A_{\xi_{\alpha}(t)}X~|_{\gamma(t)}=-[m,
\xi_{\alpha}(t)]^{\top}$. It follows from the equality $exp(tm)\cdot
[m, \xi_{\alpha}]=[m, exp(tm)\cdot \xi_{\alpha}]$ that
$A_{\xi_{\alpha}(t)}X~|_{\gamma(t)}=-exp(tm)\cdot[m, \xi_{\alpha}]^{\top}$.
Thus we obtain immediately
\begin{equation*}\label{homog shape operator}
 \nabla_X(A_{\alpha}X)~|_{z_0}=-(\frac{d}{dt}|_{t=0}[m,
\xi_{\alpha}(t)]^{\top})^{\top}=-[m, [m, \xi_{\alpha}]^{\top}]^{\top}, \quad A_{\alpha}(\nabla_XX)=-[m',
\xi_{\alpha}]^{\top}.
\end{equation*}
In this way, an equivalent condition of (2.1) for the orbit $M$ to
be an $\mathcal{A}$-manifold can be stated as
\begin{equation}\label{homog class A}
\frac{1}{2}(\nabla_X\tau)(X,X)=\sum_{\alpha=1}^p\langle~ [m,
\xi_{\alpha}]^{\top}, [m, [m, \xi_{\alpha}]^{\top}]^{\top}\rangle -\langle~ [m, \xi_{\alpha}]^{\top}, [m',
\xi_{\alpha}]^{\top}\rangle=0.
\end{equation}
\vspace{2mm}

\subsubsection{\textbf{$(m_1, m_2)=(2, 2)$}}\*
\vspace{2mm}

In this case, $G=SO(5)\times SO(5)$, $K=SO(5)$. Notice that
$$\mathfrak{k}=\Big\{\left(\begin{matrix} X&0\\
0&X \end{matrix}\right)~:~X\in\mathfrak{so}(5)\Big\}\cong\mathfrak{so}(5), \quad \mathfrak{p}=\Big\{\left(\begin{matrix} X&0\\
0&-X \end{matrix}\right)~:~X\in\mathfrak{so}(5)\Big\}\cong\mathfrak{so}(5),$$
for simplicity, we will just write the upper triangular part of a matrix in this subsection.

The group $K$ acts on $\mathfrak{p}$ by the adjoint action:
\begin{eqnarray*}
  K \times \mathfrak{p} &\rightarrow& \mathfrak{p} \\
  A,~~~Z &\mapsto& A\cdot Z\cdot A^{-1}
\end{eqnarray*}

By virtue of \cite{QTY}, the singular orbit (focal submanifold)
diffeomorphic to $\widetilde{G_2}(\mathbb{R}^5)$ is Einstein, thus
automatically Ricci parallel, and an $\mathcal{A}$-manifold.
Therefore, we will concentrate on the other singular orbit (focal
submanifold) in this subsection.

Choose a point $z_0=\frac{1}{\sqrt{2}}\left(\begin{matrix}
J&&\\
&J&\\
&&0\end{matrix}\right)$, with
$J=:\left(\begin{matrix}0&1\\-1&0\end{matrix}\right)$. It is easy to
see that the orbit $\{ A\cdot z_0\cdot A^{-1}~|~A\in SO(5)\}$
denoted by $M_1$ is a singular orbit, which is diffeomorphic to
$\mathbb{C}P^3$, as pointed out in \cite{QTY}.

A direct calculation shows that
$$\mathfrak{k}_{z_0}=\mathfrak{u}(1)\times\mathfrak{u}(1),\quad
\mathfrak{m}\cong\Big\{m=\left(\begin{matrix}
0&A&\textbf{b}\\
*&0&\textbf{c}\\
*&*&0
\end{matrix}
\right), A=\left(\begin{matrix}a_1&a_2\\a_2&-a_1\end{matrix}\right),
\textbf{b}, \textbf{c}\in M_{2,1}(\mathbb{R})
\Big\}.$$ Noticing that $JA=-AJ$ and $A^t=A$, a tangent vector of
$M_1$ at $z_0$ can be expressed as:
$$[m, z_0]=mz_0-z_0m=\frac{1}{\sqrt{2}}\left(\begin{matrix}
0&2AJ&-J\textbf{b}\\
*&0&-J\textbf{c}\\
*&*&0
\end{matrix}
\right).$$
Thus any normal vector $\xi\in T_{z_0}^{\perp}M_1\subset T_{z_0}\mathfrak{p}$ with $\langle \xi, z_0\rangle=0$ can be written as
$$\xi=\left(\begin{matrix}
\lambda J&X&0\\
*&-\lambda J&0\\
*&*&0
\end{matrix}
\right),~\mathrm{with}~ \lambda\in\mathbb{R}, JX=XJ.$$

Meanwhile, the equality
$$[m, [m, z_0]]^{\top}=\frac{1}{\sqrt{2}}\left(\begin{matrix}
0&0&3JA\textbf{c}\\
*&0&-3JA\textbf{b}\\
*&*&0
\end{matrix}
\right)=:[m', z_0]$$
implies that
$$m'=3\left(\begin{matrix}
0&0&-A\textbf{c}\\
*&0&A\textbf{b}\\
*&*&0
\end{matrix}
\right).$$
Furthermore, from
$$[m,\xi]=\left(\begin{matrix}
XA-AX^t&0&-(\lambda J\textbf{b}+X\textbf{c})\\
*&X^tA-AX&\lambda J\textbf{c}+X^t\textbf{b}\\
*&*&0
\end{matrix}
\right),
[m', \xi]=3\left(\begin{matrix}
0&0&\lambda JA\textbf{c}-XA\textbf{b}\\
*&0&\lambda JA\textbf{b}-X^tA\textbf{c}\\
*&*&0
\end{matrix}
\right),$$ and
$$[m,\xi]^{\top}=\left(\begin{matrix}
0&0&-(\lambda J\textbf{b}+X\textbf{c})\\
*&0&\lambda J\textbf{c}+X^t\textbf{b}\\
*&*&0
\end{matrix}
\right),\quad [m',\xi]^{\top}=[m', \xi],$$ it follows the equality
as below:
\begin{equation}\label{1st inner product}
  \langle~ [m, \xi]^{\top}, [m', \xi]^{\top}\rangle=0.
\end{equation}

\vspace{1mm}

On the other hand, we have
$$[m, [m, \xi]^{\top}]^{\top}=\left(
\begin{matrix}
0&*&\lambda AJ\textbf{c}+AX^t\textbf{b}\\
*&0&\lambda AJ\textbf{b}+AX\textbf{c}\\
*&*&0
\end{matrix}\right).$$
Then a simple calculation leads to
\begin{equation}\label{2nd inner product}
\langle~ [m, \xi]^{\top}, [m, [m, \xi]^{\top}]^{\top}\rangle=0.
\end{equation}

Consequently, combining (\ref{1st inner product}) with (\ref{2nd
inner product}), the proof of (\ref{homog class A}) is accomplished.
Namely, the focal submanifold diffeomorphic to $\mathbb{C}P^3$ is an
$\mathcal{A}$-manifold.

\subsubsection{\textbf{$(m_1, m_2)=(4, 5)$}}\*
\vspace{2mm}

In this case, $G=SO(10)$, $K=U(5)$, $\mathfrak{p}=\mathfrak{so}(5, \mathbb{C}).$
$K$ acts on $\mathfrak{p}$ by the adjoint action:
\begin{eqnarray*}
  K \times \mathfrak{p} &\rightarrow& \mathfrak{p} \\
  g,~~~Z &\mapsto& \overline{g}\cdot Z\cdot g^{-1}
\end{eqnarray*}

\noindent
\textbf{(1).} Choose a point $z_0=\frac{1}{\sqrt{2}}\left(\begin{matrix}
J&&\\
&J&\\
&&0\end{matrix}\right)$, with
$J=:\left(\begin{matrix}0&1\\-1&0\end{matrix}\right)$. It is easily
seen that the corresponding orbit $$M_1^{14}=\{\overline{g}\cdot
z_0\cdot g^{-1}~|~g\in U(5)\}$$ is a singular orbit (focal
submanifold), which is diffeomorphic to $U(5)\Big/Sp(2)\times U(1)$
(cf. \cite{QTY}).

Since the action of $U(5)$ on $\mathfrak{so}(5, \mathbb{C})$ is
given by $\overline{g}\cdot Z\cdot g^{-1}$, we emphasis that the
expressions before Subsection 2.2.1 for tangent space, normal space,
shape operator, and connection are still valid, only to replace the
expression of $[~,~]$ with $[m, z_0]=\overline{m}z_0-z_0m$. In this
way, the equality $exp(tm)\cdot [m, \xi_{\alpha}]=[m, exp(tm)\cdot
\xi_{\alpha}]$ still holds.

With no difficulty, we obtain that
$$\mathfrak{m}=\Big\{\left(\begin{matrix}
\lambda I&A& \textbf{b}\\
*&\mu I& \textbf{c}\\
*&*&0
\end{matrix}\right)~:~\lambda, \mu\in\sqrt{-1}\mathbb{R}, A\in gl(2, \mathbb{C}), \overline{A}J=-JA, \textbf{b},\textbf{c}\in M_{2,1}(\mathbb{C})\Big\}.$$
Notice that a tangent vector at $z_0$ can be given by
\begin{equation*}
[m, z_0]=\overline{m}z_0-z_0m=\frac{1}{\sqrt{2}}\left(\begin{matrix}
-2\lambda J&2\overline{A}J& -J\textbf{b}\\
*&-2\mu J& -J\textbf{c}\\
*&*&0
\end{matrix}\right).
\end{equation*}
Any normal vector $\xi\in T_{z_0}^{\perp}M_1\subset T_{z_0}\mathfrak{p}$ with $\langle \xi, z_0\rangle=0$ can be written as
$$\xi=\left(\begin{matrix}
tJ&X& 0\\
*&-tJ& 0\\
*&*&0
\end{matrix}\right), ~\mathrm{with}~ t\in \mathbb{R}, \overline{X}J=JX.$$
Additionally, the following equation
\begin{equation*}\label{3rd}
[m, [m, z_0]]^{\top}=\frac{1}{\sqrt{2}}\left(\begin{matrix}
0&0& 3(\lambda J\textbf{b}+JA\textbf{c})\\
*&0& 3(\mu J\textbf{c}-J\overline{A}^t\textbf{b})\\
*&*&0
\end{matrix}\right)=:[m', z_0],
\end{equation*}
implies that
$$m'=\left(\begin{matrix}
0&0& -3(\lambda \textbf{b}+A\textbf{c})\\
*&0& -3(\mu \textbf{c}-\overline{A}^t\textbf{b})\\
*&*&0
\end{matrix}\right).$$
Based on the condition of $\overline{A}J=-JA$ and
$\overline{X}J=JX$, we obtain:
\begin{equation*}\label{1}
  [m, \xi]^{\top}=\left(\begin{matrix}
-2t\lambda J+X\overline{A}^t-\overline{A}X^t&-(\lambda+\mu)X& -(tJ\textbf{b}+X\textbf{c})\\
*&2t\mu J+X^tA-A^tX& tJ\textbf{c}+X^t\textbf{b}\\
*&*&0
\end{matrix}\right)
\end{equation*}
and
\begin{equation*}\label{2}
  [m', \xi]^{\top}=3\left(\begin{matrix}
0&0& t\lambda J\textbf{b}+tJA\textbf{c}+\mu X\textbf{c}-X\overline{A}^t\textbf{b}\\
*&0& -t\mu J\textbf{c}+tJ\overline{A}^t\textbf{b}-\lambda X^t\textbf{b}-X^tA\textbf{c}\\
*&*&0
\end{matrix}\right).
\end{equation*}
Moreover, a direct calculation leads to
\begin{equation}\label{3}
\langle~ [m, \xi]^{\top} , [m', \xi]^{\top}\rangle=0.
\end{equation}
\vspace{1mm}

Next, to simplify the calculation of $\langle~ [m, \xi]^{\top},
[m, [m, \xi]^{\top}]^{\top}\rangle$, we will choose a normal
basis such that it satisfies either (i) or (ii) as follows:
\vspace{1mm}

\noindent
(i). $X=0, t=1$.
On this condition, we have
\begin{equation*}\label{4}
[m, z_0]=\frac{1}{\sqrt{2}}\left(
\begin{matrix}
-2\lambda J&2\overline{A}J&-J\textbf{b}\\
*&-2\mu J&-J\textbf{c}\\
*&*&0
\end{matrix}\right),
\quad \mathrm{and}\quad[m, \xi]^{\top}=[m, \xi]=\left(
\begin{matrix}
-2\lambda J&0&-J\textbf{b}\\
*&2\mu J&J\textbf{c}\\
*&*&0
\end{matrix}\right),
\end{equation*}
which imply
$$[m, [m, \xi]^{\top}]^{\top}=\left(
\begin{matrix}
0&*&3\lambda J\textbf{b}+\overline{A}J\textbf{c}\\
*&0&-3\mu J\textbf{c}+\overline{A}^tJ\textbf{b}\\
*&*&0
\end{matrix}\right).$$
By a simple calculation, we obtain
\begin{equation}\label{5}
\langle~ [m, \xi]^{\top}, [m, [m, \xi]^{\top}]^{\top}\rangle=0.
\end{equation}

\vspace{2mm}

\noindent
(ii). $t=0$.
On this condition, we have
\begin{equation*}
[m, z_0]=\frac{1}{\sqrt{2}}\left(
\begin{matrix}
-2\lambda J&2\overline{A}J&-J\textbf{b}\\
*&-2\mu J&-J\textbf{c}\\
*&*&0
\end{matrix}\right),
\quad[m, \xi]^{\top}=\left(
\begin{matrix}
X\overline{A}^t-\overline{A}X^t&-(\lambda+\mu)X&-X\textbf{c}\\
*&X^tA-A^tX&X^t\textbf{b}\\
*&*&0
\end{matrix}\right).
\end{equation*}
For clarity, defining
$$\sigma:=X\overline{A}^t-\overline{A}X^t,\quad \theta=:X^tA-A^tX,$$
we have
$$[m, [m, \xi]^{\top}]=\left(
\begin{matrix}
\left(
\begin{matrix}
-(3\lambda+\mu)\sigma\\
+\overline{\textbf{b}}\textbf{c}^t{X}^t-X\textbf{c}\overline{\textbf{b}}^t
\end{matrix}
\right)
&
\left(\begin{matrix}
(\lambda+\mu)^2X+\overline{A}\theta-\sigma A\\
-\overline{\textbf{b}}{\textbf{b}}^tX-X\textbf{c}\overline{\textbf{c}}^t
\end{matrix}
\right)
&
\left(\begin{matrix}
(2\lambda+\mu)X\textbf{c}\\
+(2\overline{A}X^t-X\overline{A}^t)\textbf{b}
\end{matrix}
\right)
\\
*&
\left(\begin{matrix}
-(\lambda+3\mu)\theta\\
+{X}^t\textbf{b}\overline{\textbf{c}}^t-\overline{\textbf{c}}{\textbf{b}}^tX
\end{matrix}
\right)
&
\left(\begin{matrix}
-(\lambda+2\mu)X^t\textbf{b}\\
+(2A^tX-X^tA)\textbf{c}
\end{matrix}
\right)\\
*&*&2\textbf{b}^tX\textbf{c}-2\textbf{c}^tX^t\textbf{b}
\end{matrix}\right).$$
Then a complicated but not difficult calculation shows that
\begin{equation}\label{6}
\langle~ [m, \xi]^{\top}, [m, [m, \xi]^{\top}]~\rangle=0.
\end{equation}

Finally, combining (\ref{3}) (\ref{5}) with (\ref{6}), we achieve the equality in (\ref{homog class A}).
Namely, the focal submanifold $M_1^{14}$ with $(m_1, m_2)=(4, 5)$ is an $\mathcal{A}$-manifold.

\vspace{4mm}

\noindent
\textbf{(2).} Choose a point $z_0=\left(\begin{matrix}
J&&\\
&0&\\
&&0\end{matrix}\right)$,
with $J=:\left(\begin{matrix}0&1\\-1&0\end{matrix}\right)$.
It is easily seen that the corresponding orbit
$$M_2^{13}=\{\overline{g}\cdot z_0\cdot g^{-1}~|~g\in U(5)\}$$
is a focal submanifold, which is diffeomorphic to $U(5)\Big/SU(2)\times U(3)$ (cf. \cite{QTY}).

Without much difficulty, we observe
$$\mathfrak{m}=\Big\{\left(\begin{matrix}
\lambda I&\textbf{A}\\
*&0
\end{matrix}\right)~:~\lambda\in\sqrt{-1}\mathbb{R}, \textbf{A}\in M_{2,3}(\mathbb{C})\Big\}.$$
Then a tangent vector at $z_0$ is given by:
\begin{equation*}
[m, z_0]=\overline{m}z_0-z_0m=-\left(\begin{matrix}
2\lambda J&J\textbf{A}\\
*&0
\end{matrix}\right),
\end{equation*}
and any normal vector $\xi\in T_{z_0}^{\perp}M_2^{13}\subset T_{z_0}\mathfrak{p}$ with $\langle \xi, z_0\rangle=0$ can be expressed as
$$\xi=\left(\begin{matrix}
0&0\\
0&X
\end{matrix}\right), ~\mathrm{with}~ X+X^t=0, X\in gl(3, \mathbb{C}).$$

Additionally, the equality
\begin{equation*}\label{3rd}
[m, [m, z_0]]^{\top}=\left(\begin{matrix}
0&3\lambda J\textbf{A}\\
*&0
\end{matrix}\right)=:[m', z_0],
\end{equation*}
implies that
$$m'=\left(\begin{matrix}
0&-3\lambda \textbf{A}\\
*&0
\end{matrix}\right).$$

Furthermore, we get
\begin{equation*}\label{1}
  [m, \xi]^{\top}=\left(\begin{matrix}
0&\overline{\textbf{A}}X\\
*&0
\end{matrix}\right),
\quad
\mathrm{and}
\quad [m', \xi]^{\top}=3\lambda\left(\begin{matrix}
0&\overline{\textbf{A}}X\\
*&0
\end{matrix}\right),
\end{equation*}
which leads directly to
\begin{equation}\label{7}
\langle~ [m, \xi]^{\top} , [m', \xi]^{\top}\rangle=0.
\end{equation}

On the other hand, we have
$$[m, [m, \xi]^{\top}]=\left(
\begin{matrix}
*&-\lambda \overline{\textbf{A}}X\\
*&*
\end{matrix}\right),$$
which implies immediately
\begin{equation}\label{8}
\langle~ [m, \xi]^{\top}, [m, [m, \xi]^{\top}]~\rangle=0.
\end{equation}

Finally, combining (\ref{7}) with (\ref{8}), we achieve the proof of (\ref{homog class A}),
which means that, the focal submanifold $M_2^{13}$ with $(m_1, m_2)=(4, 5)$ is an $\mathcal{A}$-manifold.


\section{\textbf{Ricci parallelism of the homogeneous cases}}

At the beginning of this section, we recall some facts for
a Riemannian manifold $M^n$ with $\pi_1M=0$.

Given $p\in M^n$, define the Ricci operator $S_p: T_pM\rightarrow
T_pM $ by $\langle S_p(X), Y\rangle=\rho(X, Y)$, $\forall~ Y\in
T_pM$. Clearly, the Ricci operator $S_p$ is a self-adjoint operator
with eigenvalues at $p$:
$$\lambda_1<\lambda_2<\cdots<\lambda_k,~~~~~ 1\leq k\leq n.$$
In this regard, we can decompose $T_pM$ into the eigenspaces $E_i$
for $S_p$ as
$$T_pM^n=E_1\oplus E_2\oplus\cdots \oplus E_k.$$

Now suppose $M^n$ is Ricci parallel, which means that the Ricci tensor is invariant under parallel
translation. Then the Ricci operator has eigenvalues
$\lambda_1<\cdots<\lambda_k$ at each point.
As a result, we can parallel translate these eigenspaces to get a
global decomposition
$$TM^n=\zeta_1\oplus \zeta_2\oplus \cdots\oplus \zeta_k,$$
into parallel distributions, with the property that
$$S|_{\zeta_i}=\lambda_i\cdot Id.$$

By the assumption $\pi_1M^n=0$, and de Rham decomposition theorem,
we can derive a global isometric splitting of $M^n$ as
$$M^n\cong N_1\times N_2\times\cdots\times N_k,~\mathrm{with} ~N_i \mathrm{~Einstein~ and}~ TN_i=\zeta_i ~(i=1, 2,\cdots, k).$$

As we mentioned in Remark 1.2, the focal submanifold $M_1$ (resp. $M_2$) with
$m_2>1$ (resp. $m_1>1$) is simply-connected.

\begin{lem}\label{k=2}
Suppose the focal submanifold $M_1^{m_1+2m_2}$ (resp. $M_2^{2m_1+m_2}$) with $g=4$ and $m_2>1$
(resp. $m_1>1$) is Ricci parallel, and the Ricci operator has eigenvalues
$\lambda_1<\cdots<\lambda_k$ with $k\geq 2$. Then $k=2$.
\end{lem}

\noindent
\emph{Proof.} We are mainly concerned with the proof for $M_1^{m_1+2m_2}$; the other case
is verbatim with obvious changes on index ranges.

Suppose $k\geq 3.$ Then a splitting for $M_1$ can be decomposed as a
product of closed manifolds: $M_1\cong N_1^{n_1}\times
N_2^{n_2}\times \overline{N}_3^{n_3}$ with $n_1\leq n_2\leq n_3$,
where $\overline{N}_3=N_3\times\cdots\times N_k$. Then from the
assumption
$0=\pi_1M_1\cong\pi_1N_2\oplus\pi_1N_2\oplus\pi_1\overline{N}_3$, we
observe that $N_1$, $N_2$ and $\overline{N}_3$ are simply-connected
as well. Thus $n_i\geq 2$ $(i=1,2,3)$.

As a matter of fact, for the focal submanifold $M_1^{m_1+2m_2}$ with
$m_2>1$, the Betti numbers satisfy (cf. \cite{Mun}):
\begin{itemize}
\item[(1)] $\beta_i(M_1)=0, ~\mathrm{for}~i\neq 0, m_2, m_1+m_2$, or $m_1+2m_2$;
\item[(2)] $\beta_j(M_1)= 1$, for $j=0, m_2, m_1+m_2$, and $m_1+2m_2$.
\end{itemize}
In fact, the homology groups of $M_1$ have no torsion. On the
condition of this fact, for $2\leq n_1\leq n_2\leq n_3\leq
m_1+2m_2-4$, from the K\"{u}nneth formula for homology group with
$\mathbb{Z}$-coefficients:
$$H_k(P\times Q)\cong\sum_{p+q=k}H_p(P)\otimes H_q(Q)\oplus\sum_{r+s=k-1}Tor(H_r(P), H_s(Q)),$$
it follows that
$$\beta_{n_1}(M_1), \beta_{n_2}(M_1), \beta_{n_3}(M_1)\geq 1, ~~\beta_{n_1+n_2}(M_1), \beta_{n_1+n_3}(M_1), \beta_{n_2+n_3}(M_1)\geq 1. $$
Consequently, we obtain that $n_1=m_2$, and thus $n_1+n_2=m_1+m_2$.
Thus $M_1^{m_1+2m_2}=N_1^{m_2}\times N_2^{m_1}\times
\overline{N}_3^{m_2}$, and further $m_1=m_2$. Moreover, it follows
that $\beta_{2m_1}(M_1)\geq 3$, a contradiction.\hfill$\Box$

\vspace{5mm} Now suppose $k=2$. From the argument above, it follows
that $M_1^{m_1+2m_2}\cong N_1^{m_2}\times N_2^{m_1+m_2}$ with
$\pi_1N_1=\pi_1N_2=0$. By K\"{u}nneth formula, we get $H_iN_1=0$ for
$1\leq i\leq m_2-1$. Similarly $H_iN_2=0$ for
$1\leq i\leq m_1-1$. In other words,  $N_1,N_2$ are all simply connected homology spheres.
Therefore, we obtain the following proposition:

\begin{prop}\label{homog RP}
Suppose the focal submanifold $M_1^{m_1+2m_2}$ with $g=4$ and $m_2>1$ is Ricci parallel, but not Einstein.
Then
\begin{itemize}
\item[(i)] The~ Ricci~ operator ~$S$~ has~ exactly~ two~ eigenvalues, ~with ~multiplicities ~$m_2$ ~and~ $m_1+m_2$, respectively;
\item[(ii)] $M_1^{m_1+2m_2}$ is diffeomorphic to a product $N^{m_2}_1 \times N^{m_1+m_2}_2$, where each factor is a simply connected homology sphere.
\end{itemize}

\end{prop}

\vspace{5mm}

In the following, we will verify the Ricci parallelism for the focal
submanifolds with $(m_1, m_2)=(2, 2)$ and $(4, 5)$ case by case. As
mentioned before, the focal submanifold diffeomorphic to
$\widetilde{G}_2(\mathbb{R}^5)$ with $(m_1, m_2)=(2, 2)$ is
Einstein, while the other one diffeomorphic to $\mathbb{C}P^3$ is
not. \vspace{3mm}

\noindent \textbf{Case 1}: The focal submanifold $M_1^6$
(diffeomorphic to $\mathbb{C}P^3$) with $(m_1, m_2)=(2, 2)$.

Suppose $M_1^6\underset{diffeo.}{\cong}\mathbb{C}P^3$ is Ricci
parallel. Then from Proposition \ref{homog RP}(ii), it follows that
$\mathbb{C}P^3$ is diffeomorphic to $N^2_1\times N^4_2$. Observe that
$N_1$ is a simply connected surface, it must be homeomorphic to a $2$-sphere, and thus
the third homotopy group $\pi_3(N_1)\cong\mathbb{Z}$. This implies that
$\pi_3(N_1^2\times N^4_2) \ncong 0$, while $\pi_3\mathbb{C}P^3 \cong 0$.
There comes a contradiction.

Consequently, the focal submanifold $M_1^6$ in this case is not Ricci parallel.

\vspace{3mm}

\noindent \textbf{Case 2}: The focal submanifold $M_1^{14}$
(diffeomorphic to $U(5)\Big/Sp(2)\times U(1)$) with $(m_1, m_2)=(4,
5)$.

Suppose that $M_1^{14}$ is Ricci parallel. It follows from
Proposition \ref{homog RP} (ii) that $M_1$ is diffeomorphic to
$N_1^5\times N_2^9$. We are going to show this impossible.

By Lemma 1.1 in \cite{Tan}, the Stiefel-Whitney class $w_4(M_2)$
of $M_2^{13}$ is nonzero (based on the elegant work of U. Abresch).
It follows that the normal bundle $\nu(M_2)$ of $M_2^{13}$ in
$S^{19}$ has $w_4(\nu(M_2))\neq 0$. By Thom isomorphism, we see
clearly that the Steenord square
$$Sq^4: H^5(M_1; \mathbb{Z}_2)\longrightarrow H^9(M_1; \mathbb{Z}_2)$$
is nonzero (compare with page 262 in \cite{Fan}), from which we claim that $M_1^{14}$ is not diffeomorphic to
$N_1^5\times N_2^9$. \vspace{3mm} To show the claim above, we
choose generators
$$e_i\in H^0(N_i;\mathbb{Z}_2), \;\;a_i \in H^{d_i}(N_i;\mathbb{Z}_2), \;i=1,2,$$
where $d_i:=\dim(N_i)$.
Denote by $p_i$ the projection from $N_1\times N_2$ to $N_i$. Then by K\"{u}nneth formula,
$p_1^*(a_1)\cup p_2^*(e_2)$ generates $H^5(M_1;\mathbb{Z}_2)$. By Cartan formula (cf. \cite{MS}), we see
$$ Sq^4(p_1^*(a_1)\cup p_2^*(e_2))=p_1^*Sq^4(a_1)\cup  p_2^*(e_2)=0 \cup p_2^*(e_2)=0, $$
a contradiction.

\noindent \textbf{Case 3}: The focal submanifold $M_2^{13}$
(diffeomorphic to $U(5)\Big/SU(2)\times U(3)$) with $(m_1, m_2)=(4,
5)$.

As in Subsection 2.2.2 (2), we choose the point $z_0=\left(\begin{matrix}
J&&\\
&0&\\
&&0\end{matrix}\right)$, with
$J=:\left(\begin{matrix}0&1\\-1&0\end{matrix}\right)$. Recall that
$M_2^{13}$ is the orbit of the isotropy representation at $z_0$.
Take a tangent vector $X=[m, z_0]=-\left(\begin{matrix}
2\lambda J&J\textbf{A}\\
*&0
\end{matrix}\right)$, with $m=\left(\begin{matrix}
\lambda I&\textbf{A}\\
*&0
\end{matrix}\right)\in \mathfrak{u}(5)$, $\lambda\in\sqrt{-1}\mathbb{R}$;
and normal vectors $\xi_{\alpha}$ with $\langle \xi_{\alpha},
z_0\rangle=0$ as $\xi_{\alpha}=\left(\begin{matrix}
0&0\\
0&X_{\alpha}
\end{matrix}\right) ~\mathrm{with}~ X_{\alpha}+X_{\alpha}^t=0, X_{\alpha}\in gl(3, \mathbb{C})$, $\alpha=1,\cdots, 6$.

Since the Ricci tensor with respect to $X, Y \in T_{z_0}M_2$ is
$\rho (X, Y) = \langle 12X-\sum_{\alpha=1}^6A_{\alpha}^2X,~
Y\rangle,$ the Ricci operator can be written as
$$\displaystyle S(X)=12X-\sum_{\alpha=1}^6A_{\alpha}^2X.$$ We are now left to complete the verification by virtue of
Proposition \ref{homog RP} (i).

From the formula $A_{\alpha}X=-[m, \xi_{\alpha}]^{\top}=:[\overline{m}_{\alpha}, z_0]$,
it follows that
\begin{equation*}
 \overline{m}_{\alpha}=\left(\begin{matrix}
0&-J\overline{\textbf{A}}X_{\alpha}\\
*&0
\end{matrix}\right)\in \mathfrak{u}(5).
\end{equation*}
Thus choosing
$X_1=\left(\begin{smallmatrix}
0&1&0\\
-1&0&0\\
0&0&0
\end{smallmatrix}\right)$,
$X_2=\left(\begin{smallmatrix}
0&i&0\\
-i&0&0\\
0&0&0
\end{smallmatrix}\right)$,
$X_3=\left(\begin{smallmatrix}
0&0&1\\
0&0&0\\
-1&0&0
\end{smallmatrix}\right)$,
$X_4=\left(\begin{smallmatrix}
0&0&i\\
0&0&0\\
-i&0&0
\end{smallmatrix}\right)$,
$X_5=\left(\begin{smallmatrix}
0&0&0\\
0&0&1\\
0&-1&0
\end{smallmatrix}\right)$, and
$X_6=\left(\begin{smallmatrix}
0&0&0\\
0&0&i\\
0&-i&0
\end{smallmatrix}\right)$, we derive that
\begin{equation*}
\sum_{\alpha=1}^6A_{\alpha}^2X = -\sum_{\alpha=1}^6[\overline{m}_{\alpha}, \xi_{\alpha}]^{\top}
= \sum_{\alpha=1}^6\left(\begin{matrix}
0&-J\textbf{A}\overline{X}_{\alpha}X_{\alpha}\\
*&0
\end{matrix}\right)^{\top}
= \left(\begin{matrix}
0&4J\textbf{A}\\
*&0
\end{matrix}\right).
\end{equation*}
In this way, we obtain the Ricci operator
$$S(X)=12X-\sum_{\alpha=1}^6A_{\alpha}^2X=\left(\begin{matrix}
-24\lambda J&-8J\textbf{A}\\
*&0
\end{matrix}\right).$$
A direct calculation shows that the Ricci operator $S$ has two
eigenvalues $12$ and $8$, with multiplicities $1$ and $12$,
respectively, which contradicts Proposition \ref{homog RP} (i).

The proof of Theorem \ref{Ricci parallel} (ii), (iii) is now
complete.


\vspace{4mm}

\section{\textbf{Ricci parallelism of OT-FKM type}}
For convenience, we will firstly deal with the focal submanifold $M_2$ of OT-FKM type.

\subsection{$M_2$ of OT-FKM type}

The proof of Theorem \ref{Ricci parallel} will be finished by establishing the following two propositions.

\begin{prop}\label{M2 m=1}
The focal submanifold $M_2$ of OT-FKM type with $m=1$ is Ricci parallel.
\end{prop}

\noindent
\emph{Proof.} When $m=1$, the equalities (\ref{AX AY}) turn to be
\begin{equation}\label{AX AY m=1}
\tau(X, Y)=\sum_{\alpha=1}^{l-1}\langle A_{\alpha}X,
A_{\alpha}Y\rangle=\langle X, Y\rangle+(l-2)V(X, Y),
\end{equation}
where $V(X, Y)$ is defined in (\ref{V W}). At $x\in M_2$ with
$Px=x$, we can always choose $Q_0=:P=\langle P_0x, x\rangle P_0 +
\langle P_1x, x\rangle P_1$ by the definition of FKM-polynomial $F$
in (\ref{FKM isop. poly.}), and then $Q_1$ can be stated as
$Q_1=\langle P_1x, x\rangle P_0 - \langle P_0x, x\rangle P_1$. In
this way, it is easily seen that $$\displaystyle V(X,
Y)=\sum_{k=0}^1 \langle X, P_kx\rangle\langle Y, P_kx\rangle=\langle
X, Q_1x\rangle\langle Y, Q_1x\rangle,$$ and then
\begin{equation}\label{m=1 M2 V}
\nabla V=0\Longleftrightarrow \langle X, \nabla_Z(Q_1x)\rangle\langle X, Q_1x\rangle =0, ~~\forall X, Z\in T_xM_2.
\end{equation}
For the first factor on the righthand side, it follows from a simple
calculation that
\begin{eqnarray*}
\langle X, \nabla_Z(Q_1x)\rangle &=& \langle X,~ D_Z(\langle P_1x, x\rangle P_0x - \langle P_0x, x\rangle P_1x)\rangle \\
   &=& \langle X, ~Q_1Z+2(\langle P_1x, Z\rangle P_0x - \langle P_0x, Z\rangle P_1x)\rangle \\
   &=& 0.
\end{eqnarray*}
We need to explain the reason for the last equality. In this case
$\{Q_1N_1, Q_1N_2, \cdots, Q_1N_{l-1}, Q_1x\}$ constitutes an
orthonormal basis of $T_xM_2$ by (\ref{o.n.b of TxM2}), and we can
show that
\begin{equation*}\label{Qx Ax}
 \langle Q_1Z, Q_1x\rangle=0,\quad \langle~ \langle P_1x, Z\rangle P_0x - \langle P_0x, Z\rangle P_1x, ~Q_1x~\rangle=\langle Z, Q_0x\rangle=0,
\end{equation*}
and for any $\alpha=1,\cdots, l-1$,
\begin{eqnarray*}\label{QN Ax}
&&\langle Q_1Z,
Q_1N_{\alpha}\rangle=0,\\
&&\langle~ \langle P_1x, Z\rangle P_0x - \langle P_0x, Z\rangle
P_1x, ~Q_1N_{\alpha}~\rangle=- \langle N_{\alpha},
Q_1x\rangle\langle Z, Q_1x\rangle=0.
\end{eqnarray*}

The proof of Proposition \ref{M2 m=1} is now complete. \vspace{3mm}

\begin{rem}
In fact, as asserted by \cite{TY2}, up to a two-fold covering, $M_2$
with $m=1$ is isometric to $S^1\times S^{l-1}$.
\end{rem}

\begin{prop}\label{M2 m>2}
The focal submanifold $M_2$ of OT-FKM type with $m\geq 2$ is not Ricci parallel.
\end{prop}

\noindent
\emph{Proof.} Recall the equalities (\ref{AX AY})
\begin{equation*}\label{AX AY m=2}
\tau(X, Y)=\sum_{\alpha=1}^{l-m}\langle A_{\alpha}X,
A_{\alpha}Y\rangle=m\langle X, Y\rangle+(l-2m)V(X, Y)-W(X, Y).
\end{equation*}
For covariant derivative of the items on the righthand side, we have
\begin{eqnarray}\label{m=2 M2 2nd}
  (\nabla_ZV)(X, Y)&=&\sum_{i=0}^m\langle X, P_iZ\rangle\langle Y, P_ix\rangle+ \langle X, P_ix\rangle\langle Y, P_iZ\rangle\nonumber\\
  &=& \sum_{i=1}^m\langle X, Q_iZ\rangle\langle Y, Q_ix\rangle+ \langle X, Q_ix\rangle\langle Y, Q_iZ\rangle.
\end{eqnarray}
and by (\ref{W}), it is not difficult to see
\begin{eqnarray}\label{nabla W}
&& (\nabla_ZW)(X, Y) \\
&=& \sum_{k,s=0, k\neq s}^m\Big(\langle X, P_kP_sZ\rangle\langle Y, P_kP_sx\rangle + \langle X, P_kP_sx\rangle\langle Y, P_kP_sZ\rangle\Big) - 2(\nabla_ZV)(X, Y) \nonumber\\
&=& \sum_{i,j=1, i\neq j}^m\Big(\langle X, Q_iQ_jZ\rangle\langle Y,
Q_iQ_jx\rangle + \langle X, Q_iQ_jx\rangle\langle Y,
Q_iQ_jZ\rangle\Big) \nonumber\\
&& +2\sum_{i=1}^m \Big(\langle X, Q_ix\rangle\langle Y,
Q_iQ_0Z\rangle+\langle Y, Q_ix\rangle\langle X, Q_iQ_0Z\rangle\Big)-
2(\nabla_ZV)(X, Y)\nonumber
\end{eqnarray}
Taking $X=Q_1Q_2x, Y=Q_1x$ and $Z=Q_2x$ in (\ref{m=2 M2 2nd}), (\ref{nabla W}), we obtain
\begin{equation}\label{nabla u}
  (\nabla_Z\tau)(X, Y)=l-2m+2+\sum_{\begin{smallmatrix}i,j=1,\cdots, m, i\neq j,\\ \{i, j\}\neq \{1, 2\}\end{smallmatrix}}\langle Q_1Q_2x, Q_iQ_jx\rangle^2\geq l-2m+2.
\end{equation}
Suppose $M_2$ is Ricci parallel. Then we get $l-2m+2\leq 0$, which holds only in the cases $(m_1, m_2)=(6, 1)$, $(5, 2)$ and $(9, 6)$ in OT-FKM type.

While in view of \cite{FKM}, the families with multiplicities $(6, 1)$ and $(5, 2)$ are congruent to that with multiplicities $(1, 6)$ and $(2, 5)$ respectively, and the
focal submanifolds interchange. Thus the focal submanifold $M_2$ with $(m_1, m_2)=(6, 1)$ and $(m_1, m_2)=(5, 2)$ are congruent to $M_1$ with
$(m_1, m_2)=(1, 6)$ and $(2, 5)$ respectively, which are not Ricci parallel as a direct result of (\ref{dim V}) in Subsection 4.2.

As for the $(9, 6)$ case, we have
\begin{lem}\label{(9,6)}
The focal submanifold $M_2$ of OT-FKM type with $(m_1, m_2)=(9,6)$ is not Ricci parallel.
\end{lem}
Up to now, the proof of Proposition \ref{M2 m>2} is complete provided we give a proof of the lemma above.

\noindent
\emph{Proof.}
Suppose that $M_2^{24}$ is Ricci parallel. It follows from
Proposition \ref{homog RP}(ii) that $M_2$ is diffeomorphic to
$N_1^9\times N_2^{15}$, where $N_1$ and $N_2$ are simply connected homology spheres.
By the well-known Hurewicz theorem and Whitehead theorem, one sees that any simply connected homology sphere is in fact homotopy equivalent to a unit sphere.
Thus $M_2^{24}$ has the same homotopy type with $S^9\times S^{15}$.

On the other hand, it is impossible that $M_2^{24}$ has the same homotopy type with $S^9\times S^{15}$ by the Clifford construction.
To show this claim, we follow Wang (\cite{Wan}). For a symmetric Clifford system
$\{P_0,\cdots,P_m\}$ on $\mathbb{R}^{2l}$ with $l=k\delta(m)$, where $k$ is a positive integer and $\delta(m)$ is
the dimension of irreducible $C_{m-1}$-modules, by using the theory of Atiyah-Bott-Shapiro, Wang constructed a vector bundle $\xi$ of rank $l$ over $S^m$.
Moreover, Wang showed (Prop.1 in \cite{Wan}) that the focal manifold
$M_2$ is diffeomorphic to $S(\xi)$, the unit sphere bundle of $\xi$. Suppose that $l \geq m+2$. Thus the vector $\xi$ is trivial if and only if it is stable trivial.
If $m$ is not divisible by $4$, observe (ref. the proof of Corollary $1$ in \cite{Wan}) that $\xi-l\in \widetilde{KO}(S^m)$ is equal to $k$ times a generator of $\widetilde{KO}(S^m)$.

In our case, $l=16,m=9$.  The assumption $l \geq m+2$ is satisfied. Furthermore, since $\delta(9)=16$, $k=1$, it follows from the arguments above that
$$ \xi-16\in \widetilde{KO}(S^9) \cong \mathrm{Z}_2 \;\;\mathrm{is} \;\mathrm{a} \;\mathrm{generator}. $$
Thus the characteristic map $\chi(\xi)$ of the bundle $\xi$ over $S^9$ is not trivial in $\pi_8SO(16)$. Consider the $J$-homomorphism of Whitehead
$$ J: \pi_8SO(16)\cong \pi_8SO \cong \mathrm{Z}_2 \longrightarrow  \pi_{24}S^{16} \cong \pi_8^S.  $$
By Adams\cite{Adm}, the homomorphism $J$ is a monomorphism. Hence $J(\chi(\xi))$ does not vanish in  the stable homotopy group $\pi_8^S$.
Applying Theorem (1.11) in \cite{JW}, we conclude that $M_2^{24}$ has not the same homotopy type with $S^9\times S^{15}$.

\vspace{4mm}

\subsection{$M_1$ of OT-FKM type}

From (\ref{M1 class A}), it follows that the focal submanifold $M_1$ of OT-FKM type
is Ricci parallel, if and only if at any point $x\in M_1$,
\begin{equation}\label{9}
\sum_{\alpha, \beta =0, \alpha \neq \beta}^{m}\langle Y, ~~\langle
X, P_{\alpha}P_{\beta}x\rangle P_{\alpha}P_{\beta}Z+\langle Z,
P_{\alpha}P_{\beta}x\rangle P_{\alpha}P_{\beta}X\rangle=0, \quad
\forall ~X, Y, Z\in T_xM_1.
\end{equation}
It is easily seen that an equivalent condition of (\ref{9}) can be stated as:
\begin{equation*}\label{M1 of OT-FKM Ricci parallel}
\begin{array}{l}
\displaystyle A(X, Y)=:\sum_{\alpha, \beta =0, \alpha< \beta}^{m}\langle X, P_{\alpha}P_{\beta}x\rangle P_{\alpha}P_{\beta}Y+\langle X, P_{\alpha}P_{\beta}Y\rangle P_{\alpha}P_{\beta}x\\
\displaystyle \qquad\qquad\in \mathbb{R}x\oplus Span\{P_0x, P_1x,\cdots, P_mx\}=:\mathcal{L},\qquad\qquad \forall~X, Y\in T_xM_1,~\forall~ x\in M_1.
\end{array}
\end{equation*}
Define $$\mathcal{V}_x=:Span\{P_{\alpha}P_{\beta}x, 0\leq\alpha<
\beta\leq m\}\subset T_xM_1, \quad
\mathcal{W}_x=:\mathcal{V}_x^{\perp}\subset T_xM_1,$$ so that
$T_xM_1=\mathcal{V}_x\oplus \mathcal{W}_x$.

Suppose now that $M_1$ is Ricci parallel. Firstly, for any $X=w\in
\mathcal{W}_x$, we have
\begin{equation}\label{A(w, Y)}
A(w, Y)=\sum_{\alpha, \beta =0, \alpha< \beta}^{m}\langle w, P_{\alpha}P_{\beta}Y\rangle P_{\alpha}P_{\beta}x\in \mathcal{L}\cap \mathcal{V}_x=\{0\}, \quad\forall ~Y\in T_xM_1.
\end{equation}
Next, choosing $Y=P_0w\in T_xM_1$, (\ref{A(w, Y)}) changes to
\begin{equation}\label{A(w,P0w)}
A(w, P_0w)=\sum_{\alpha, \beta =0, \alpha< \beta}^{m}\langle w, P_{\alpha}P_{\beta}P_0w\rangle P_{\alpha}P_{\beta}x=\sum_{\beta=1}^m\langle w, P_{\beta}w\rangle P_0P_{\beta}x=0.
\end{equation}
Since $P_0P_1x, P_0P_2x, \cdots, P_0P_mx$ are linearly independent, (\ref{A(w,P0w)}) implies that
\begin{equation*}\label{wPw=0}
\langle w, P_{\beta}w\rangle =0, \quad \beta=1,2,\cdots, m.
\end{equation*}
Analogously, replacing $Y=P_0w$ with $Y=P_1w, P_2w, \cdots, P_mw$
leads to
\begin{equation*}\label{wPaw=0}
\langle w, P_{\alpha}w\rangle =0, \quad \alpha=0, 1, \cdots, m.
\end{equation*}
Using a polarization, it is easy to find that $\langle w_1,
P_{\alpha}w_2\rangle =0$, for any $w_1, w_2\in \mathcal{W}_x$. In
other words,
$$P_{\alpha}w\in \mathcal{V}_x, \quad \mathrm{for~any}~ w\in \mathcal{W}_x.$$
Denote the shape operator with respect to $P_0x$ by $A_0=:A_{P_0x}$, then
\begin{eqnarray*}
A_0: T_xM_1 &\rightarrow& T_xM_1=T_0\oplus T_1\oplus T_{-1} \\
  X &\mapsto& -(P_0X)^{\top} \nonumber
\end{eqnarray*}
where $T_0, T_1, T_{-1}$ are eigenspaces of $A_0$ corresponding to eigenvalues $0, 1, -1$, respectively,
and in this case, $T_0=Span \{P_0P_1x, P_0P_2x, \cdots, P_0P_mx\}\subset \mathcal{V}_x$.
Thus
\begin{eqnarray*}\label{A0 W}
A_0~|_{\mathcal{W}_x}: \mathcal{W}_x &\rightarrow& T_0^{\perp}\subset \mathcal{V}_x\\
  w&\mapsto& A_0w=-P_0w\in \mathcal{V}_x \nonumber
\end{eqnarray*}
is injective. Then it is clear that $\dim \mathcal{W}_x\leq \dim
T_0^{\perp}$, which implies immediately the following necessary
condition for $M_1$ to be Ricci parallel:
\begin{equation}\label{dim V}
\dim \mathcal{V}_x\geq l-1=k\delta(m)-1.
\end{equation}


On the other hand, $\dim \mathcal{V}_x\leq \frac{1}{2}m(m+1)$.
Comparing with the following inequalities
$$l-1>\frac{1}{2}m(m+1), \quad l-m-1>0,$$ we are left to deal with  the following cases, while the others are not Ricci parallel:
\vspace{2mm}

\begin{table}[htbp]
\begin{tabular}{|c|c|c|c|c|c|c|c|c|c|c|}
  \hline
  $m$ & 2 & 4 & 5 & 6 & 7 & 8 & 9 & 10 & 11 & 12 \\\hline
  $k$ & 2 & 2 & 1,2 & 1,2 & 2,3 & 2,3,4 & 1,2 & 1 & 1 & 1 \\
  \hline
\end{tabular}
\end{table}

\noindent \textbf{(1) the case $m=2, k=2$, i.e. $(m_1, m_2)=(2,
1)$.}

In view of \cite{FKM}, the family with multiplicities $(2, 1)$ is congruent to that with multiplicities $(1, 2)$, and the
focal submanifolds interchange. Thus the focal submanifold $M_1$ with $(m_1, m_2)=(2, 1)$ is congruent to $M_2$ with
$(m_1, m_2)=(1, 2)$, which is Ricci parallel according to Proposition \ref{M2 m=1}.
\vspace{2mm}

\noindent \textbf{(2) the case $m=4, k=2$, i.e. $(m_1, m_2)=(4,
3)$.}

According to \cite{FKM}, there are two examples of OT-FKM type isoparametric
polynomials with multiplicities $(m_1, m_2)=(4, 3)$, which are distinguished by an invariant
$$\mathrm{Trace}(P_0P_1P_2P_3P_4)=2q\delta(4),~\mathrm{with}~q\equiv 2 ~mod~2.$$
When $q=2$, \cite{QTY} asserts that the $M_1$ is Einstein. Thus we are left to the other case $q=0$,
in which $P_0P_1P_2P_3P_4\neq \pm Id$.

Setting $P=P_0P_1P_2P_3$, it is easy to see that $P$ is symmetric
and $P^2=Id$. Then following from Theorem 5.1 and 5.2 in \cite{FKM},
we can find a point $x\in M_1$ as the $+1$-eigenvector of $P$, i.e.
$P_0P_1P_2P_3x=x$. On this condition, we can show
$$\mathcal{V}_x=Span\{P_0P_1x, P_0P_2x, P_0P_3x, P_0P_4x, P_1P_4x, P_2P_4x, P_3P_4x\}.$$
Then from the decomposition $T_xM_1=\mathcal{V}_x\oplus \mathcal{W}_x$, it follows that
$$\mathcal{W}_x=Span\{ P_0P_1P_4x, P_0P_2P_4x, P_0P_3P_4x\}.$$

On the other hand, using polarization, another equivalent condition
of (\ref{9}) can be stated as well:
\begin{equation}\label{10}
B(X)=:\sum_{\alpha, \beta =0, \alpha< \beta}^{m}\langle X,
P_{\alpha}P_{\beta}x\rangle P_{\alpha}P_{\beta}X\in \mathcal{L},
\quad \mathrm{for~any}~X\in T_xM_1, x\in M_1.
\end{equation}
Choosing now $X=P_0P_1x+w\in \mathcal{V}_x\oplus \mathcal{W}_x$, we
get
$$B(X)=-2x+(P_0P_1-P_2P_3)w.$$

Suppose that $M_1$ is Ricci parallel. 
Noticing $\langle (P_0P_1-P_2P_3)w, x\rangle=\langle
(P_0P_1-P_2P_3)w, P_ix\rangle=0$ $(i=0,1,2,3)$, the arguments above
imply
$$(P_0P_1-P_2P_3)w  \sslash P_4x.$$
However, setting $w=P_0P_2P_4x$, we have $\langle (P_0P_1-P_2P_3)w,
P_4x\rangle=-2\langle P_1P_2P_4x, P_4x\rangle=0$. Then it must be
true that $P_1P_2P_4x=0$, a contradiction.

Therefore, $M_1$ with $(m_1, m_2)=(4, 3)$ and $P_0P_1P_2P_3P_4\neq \pm Id$ is not Ricci parallel.

\vspace{2mm}

\noindent
\textbf{(3) the cases $m=5$, $k=1, 2$, i.e. $(m_1, m_2)=(5, 2), (5, 10)$.}

Choose $x\in S^{2l-1}$ as a common eigenvector of the commuting 4-products
$P_0P_1P_2P_3$ and $P_0P_1P_4P_5$. It is easy to see that $x\in M_1$ and
$\dim \mathcal{V}_x=7$.

In the case $(m_1, m_2)=(5, 2)$, $\dim \mathcal{W}_x=2$,
and $\mathcal{W}_x=Span\{P_0P_2P_4x, P_0P_2P_5x\}$.
Suppose $M_1$ is Ricci parallel. Then for $w=P_0P_2P_4x$, $X=P_0P_2P_5x$,
we have
$$A(w, X)=\sum_{\alpha, \beta=0, \alpha<\beta}^5\langle w, P_{\alpha}P_{\beta}X\rangle P_{\alpha}P_{\beta}x=3P_4P_5x\neq 0,$$
which contradicts (\ref{A(w, Y)}).

In the case $(m_1, m_2)=(5, 10)$, $\dim \mathcal{V}_x=7<l-1=15$, which means that $M_1$ is not Ricci parallel by (\ref{dim V}).

\vspace{2mm}

\noindent
\textbf{(4) the cases $m=6$, $k=1, 2$, i.e. $(m_1, m_2)=(6, 1), (6, 9)$.}

In the case $(m_1, m_2)=(6, 1)$, according to \cite{FKM}, $M_1$ is
congruent to $M_2$ with $(m_1, m_2)=(1, 6)$, which is Ricci
parallel by Proposition \ref{M2 m=1}.

In the case $(m_1, m_2)=(6, 9)$, choose $x\in S^{2l-1}$ as a common
eigenvector of the commuting 4-products $P_0P_1P_2P_3$,
$P_0P_1P_4P_5$ and $P_0P_2P_4P_6$. Then it is easy to see that $x\in
M_1$ and $\dim \mathcal{V}_x\leq7<l-1=15$. It follows from (\ref{dim
V}) that $M_1$ is not Ricci parallel.

\vspace{2mm}

\noindent
\textbf{(5) the cases $m=7$, $k=2, 3$, i.e. $(m_1, m_2)=(7, 8), (7, 16)$.}

Choose $x\in S^{2l-1}$ as a common eigenvector of the commuting $4$-products
$P_0P_1P_2P_3$, $P_0P_1P_4P_5$, $P_0P_1P_6P_7$ and $P_0P_2P_4P_6$.
Then it is easily seen that $x\in M_1$ and $\dim \mathcal{V}_x=7$.

In these two cases, we have $l=k\delta(7)=16$ or $24$. It follows immediately that
$\dim \mathcal{V}_x<l-1$, thus $M_1$ in both cases are not Ricci parallel.

\vspace{2mm}

\noindent
\textbf{(6) the cases $m=8$, $k=2, 3, 4$, i.e. $(m_1, m_2)=(8, 7), (8, 15), (8, 23)$.}

When $k=2~ (\mathrm{resp}.~ 3, 4)$, the FKM-polynomial is defined on $\mathbb{R}^{32}$ $(\mathrm{resp}. ~\mathbb{R}^{48}, \mathbb{R}^{64})$.
Since $P_2P_4P_6P_8$ anti-commutes with $P_2$, $E_+(P_2P_4P_6P_8)$
has dimension $16$ $(\mathrm{resp}.~ 24, 32)$.
It is an invariant subspace of the anti-commuting operators $P_3P_4P_7P_8$ and $P_3$.
Thus $E_+(P_2P_4P_6P_8)\cap E_+(P_3P_4P_7P_8)$ is of dimension $8$ $(\mathrm{resp}.~ 12, 16)$
and further it is an invariant subspace of the anti-commuting operators $P_5P_6P_7P_8$ and $P_5$.
Thus the space $E'=:E_+(P_2P_4P_6P_8)\cap E_+(P_3P_4P_7P_8)\cap E_+(P_5P_6P_7P_8)$ is of dimension
$4$ $(\mathrm{resp}. ~6, 8)$ and on this space, we have
$$F(x)=|x|^4-2\sum_{\alpha=0}^1\langle P_{\alpha}x, x\rangle^2.$$

This function is not constant and a maximum point lies in $M_1$. We
choose such an $x\in M_1$. Then at this point, it is not difficult
to see that $\dim \mathcal{V}_{x}\leq 22.$

In the case $k=2~ (\mathrm{resp}. ~3)$, $l-1=23 ~(\mathrm{resp}.
~31)$, we have $\dim \mathcal{V}_{x}<l-1$. A similar argument as
above shows that $M_1$ with $(m_1, m_2)=(8, 15)$ or $(8, 23)$ is not
Ricci parallel.

In the case $k=1$, we divide the proof into two cases: the definite
family $P_0P_1\cdots P_8=\pm Id$ and the indefinite family
$P_0P_1\cdots P_8\neq\pm Id$.

Case 1: For the definite family, we observe that
$$\{P_0P_1x, \cdots, P_0P_8x,P_1P_2x,\cdots, P_1P_8x,P_2P_3x, \cdots, P_2P_8x, P_3P_4x\}$$
constitutes an orthonormal basis of $T_{x}M_1$. Taking $X=P_0P_3x$,
$Y=P_0P_2x$, we see
\begin{eqnarray*}
A(X, Y) &=&  \sum_{\alpha, \beta =0, \alpha< \beta}^{m}\langle X, P_{\alpha}P_{\beta}x\rangle P_{\alpha}P_{\beta}Y+\langle X, P_{\alpha}P_{\beta}Y\rangle P_{\alpha}P_{\beta}x\\
 &=& \sum_{\alpha, \beta =0, \alpha< \beta}^{m} \langle P_0P_3x, P_{\alpha}P_{\beta}x\rangle P_{\alpha}P_{\beta}P_0P_2x+\langle P_0P_3x, P_{\alpha}P_{\beta}P_0P_2x\rangle P_{\alpha}P_{\beta}x\\
 &=& 2P_2P_3x \\
 &\not \in& \mathcal{L}
\end{eqnarray*}
Thus $M_1$ in this case is not Ricci parallel.

Case 2: For the indefinite family, extend $\{P_0, P_1,\cdots, P_8\}$ to $\{P_0, P_1,\cdots, P_9\}$.
Choose $x$ to be a common eigenvector of $P_{2\alpha}P_{2\alpha+1}P_{2\beta}P_{2\beta+1}$, $0\leq \alpha<\beta\leq 4$.
Then $x\in M_1$ and $\dim \mathcal{V}_x=21$.
On the other hand, since $M_1$ is of dimension $22$, the Ricci operator
$S(X)=2(l-m-2)X+2\sum_{0\leq \alpha<\beta\leq 9}\langle X, P_{\alpha}P_{\beta}x\rangle P_{\alpha}P_{\beta}x$,
must have an eigenvalue $0$ with multiplicity $1$.

Suppose that $M_1$ is Ricci parallel, which is indeed not Einstein
(cf. \cite{QTY}). Then Proposition \ref{homog RP} (i) shows that the
Ricci operator must have two eigenvalues with multiplicities $7$ and
$15$, respectively. There comes a contradiction.

Therefore, $M_1$ in this case is not Ricci parallel.

\vspace{2mm}

\noindent
\textbf{(7) the cases $m=9$, $k=1, 2$, i.e. $(m_1, m_2)=(9, 6), (9, 22)$.}

In the case $(m_1, m_2)=(9, 6)$, the Ricci curvature with respect to $X, Y\in T_xM_1$
can be stated as (cf. \cite{QTY}):
\begin{equation}\label{9,6 M1 Ricci}
\rho(X, Y) = 10\langle X, Y\rangle+4\big\{\frac{5}{2}\langle X,
P_0P_1x\rangle\langle Y, P_0P_1x\rangle+\sum_{(\alpha,
\beta)\in\Lambda}\langle X, P_{\alpha}P_{\beta}x\rangle\langle Y,
P_{\alpha}P_{\beta}x\rangle\big\}
\end{equation}
where $\Lambda=\{(0, 2), (0, 3),\cdots,(0, 9), (2, 4),(2,5), \cdots, (2, 9), (4, 6), (4, 7),\cdots, (4, 9), (6, 8),
(6, 9)\}$.
By a direct calculation, we derive that
\begin{eqnarray*}
\frac{1}{4}(\nabla_Z\rho)(X, Y) &=& \frac{5}{2}\langle X, P_0P_1Z\rangle\langle Y, P_0P_1x\rangle+\frac{5}{2}\langle X, P_0P_1x\rangle\langle Y, P_0P_1Z\rangle\\
  &&  +\sum_{(\alpha, \beta)\in\Lambda}(\langle X, P_{\alpha}P_{\beta}Z\rangle\langle Y, P_{\alpha}P_{\beta}x\rangle+\langle X, P_{\alpha}P_{\beta}x\rangle\langle Y,
  P_{\alpha}P_{\beta}Z\rangle).
\end{eqnarray*}
Taking now tangent vectors $X=P_0P_1x$, $Y=P_0P_2x$ and $Z=P_1P_2x$,
we get
$$\frac{1}{4}(\nabla_Z\rho)(X, Y)=\frac{3}{2}+\sum_{\alpha=6}^9\langle P_0P_2x, P_4P_{\alpha}x\rangle^2+\sum_{\beta=8}^9\langle P_0P_2x, P_6P_{\beta}x\rangle^2\geq \frac{3}{2}.$$
Thus the $M_1$ with $(m_1, m_2)=(9, 6)$ is not Ricci parallel.

In the case $(m_1, m_2)=(9, 22)$, choose $x\in S^{2l-1}$ as a common
eigenvector of the commuting $4$-products
$P_{2\alpha}P_{2\alpha+1}P_{2\beta}P_{2\beta+1}$, $0\leq
\alpha<\beta\leq 4$. Then it is easy to see that $x\in M_1$ and
$\dim \mathcal{V}_x\leq 21$. Evidently, $\dim \mathcal{V}_x<l-1=31$,
thus $M_1$ is not Ricci parallel.

\vspace{2mm}

\noindent
\textbf{(8) the case  $m=10$, $k=1$, i.e. $(m_1, m_2)=(10, 21)$.}

With a similar discussion as in the case (6), it follows that the
space $E_+(P_0P_1P_2P_3)\cap E_+(P_0P_1P_4P_5)\cap E_+(P_4P_5P_6P_7)
\cap E_+(P_2P_3P_8P_9)$ is of dimension $4$. On this space, the
FKM-polynomial is
$$F(x)=|x|^4-2\langle P_{10}x, x\rangle^2.$$
This function is not constant and a maximum point lies in $M_1$. We choose
$x\in S^{63}$ to be the maximum point of the restricted $F$.
It is easily seen that $x\in M_1$, and $\dim \mathcal{V}_x\leq 31=l-1$.

If $\dim \mathcal{V}_x<31$, then $M_1$ is not Ricci parallel.

If $\dim \mathcal{V}_x=31=l-1$, since $M_1$ is of dimension $52$,
$0$ must be an eigenvalue of the Ricci operator $S$
with multiplicity $21$. Suppose that $M_1$ is Ricci parallel, which
is indeed not Einstein by \cite{QTY}. Then the Ricci operator $S$
has two eigenvalues with multiplicities $21$ and $31$, respectively.
Thus for any tangent vector $X\in \mathcal{V}_x$, $S(X)=c\cdot X$,
where $c$ is a constant. However, taking $X_1=P_0P_1x$, we have
$S(X_1)=5P_0P_1x$; while taking $X_2=P_0P_{10}x$, we have
$S(X_2)=P_0P_{10}x$, which is an obvious contradiction.

Therefore, $M_1$ is not Ricci parallel.

\vspace{2mm}

\noindent
\textbf{(9) the case $m=11$, $k=1$, i.e. $(m_1, m_2)=(11, 52)$.}

Choose $x\in S^{127}$ as a common eigenvector of the commuting 4-
products $P_{2\alpha}P_{2\alpha+1}P_{2\beta}P_{2\beta+1}$, $0\leq \alpha<\beta\leq 5$.
Then it is easy to see that $x\in M_1$ and further $\dim \mathcal{V}_x\leq 31<l-1=63$.
Thus $M_1$ is not Ricci parallel.

\vspace{2mm}

\noindent
\textbf{(10) the case $m=12$, $k=1$, i.e. $(m_1, m_2)=(12, 51)$.}

Choose $x\in S^{127}$ as a common eigenvector of the commuting 4-
products $P_0P_1P_2P_3$, $P_4P_5P_6P_7$, $P_0P_1P_8P_9$, $P_2P_3P_8P_9$, $P_6P_7P_{10}P_{11}$, and $P_0P_2P_{8}P_{12}$.
Then it is easy to see that $x\in M_1$ and further $\dim \mathcal{V}_x\leq 56<l-1=63$.
Thus $M_1$ is not Ricci parallel.

The proof of Theorem \ref{Ricci parallel} (i) is now complete.

\section{\textbf{Examples to the problem of Besse}}
We begin this section with a proof of Proposition \ref{Rie product}.

\noindent
\textbf{Proposition 1.1.} The focal submanifolds of isoparametric hypersurfaces in spheres with $g=4$
and $m_1, m_2>1$ are not Riemannian products.
\vspace{1mm}

\noindent
\emph{Proof.}  For convenience, we are only concerned with the proof for $M_1$, while the other case is verbatim.

Observe that the sectional curvature is given by
\begin{eqnarray*}                                                          
Sec(X\wedge Y) &=& 1+\sum_{\alpha=1}^{m_1+1}\langle A_{\alpha}X,
X\rangle\langle A_{\alpha}Y, Y\rangle -
\sum_{\alpha=1}^{m_1+1}\langle A_{\alpha}X, Y\rangle^2,
\end{eqnarray*}
where $X$ and $Y$ are unit tangent vectors at the same point
perpendicular to each other. For simplicity, we denote
$$\widetilde{A}(X, Y)=\sum_{\alpha=1}^{m_1+1}\langle A_{\alpha}X,
X\rangle\langle A_{\alpha}Y, Y\rangle,\quad \widetilde{B}(X,
Y)=\sum_{\alpha=1}^{m_1+1}\langle A_{\alpha}X, Y\rangle^2.$$
\begin{lem}\label{lem Rie prod}
The inequality $\widetilde{A}\leq 1$ holds, and the equality holds
if and only if the following two conditions are both satisfied:
\begin{itemize}
\item[(1)]$\langle A_{\alpha}X, X\rangle=\langle A_{\alpha}Y, Y\rangle$, for any $\alpha=1,\cdots,m_1+1$
\item[(2)] $X$ is an $+1$-eigenvector for a certain unit normal vector $N$.
\end{itemize}
\end{lem}

\noindent
Proof of Lemma \ref{lem Rie prod}. Denote
\begin{eqnarray*}
\textbf{a} &=:& (\langle A_1X, X\rangle, \langle A_2X, X\rangle,\cdots, \langle A_{m_1+1}X, X\rangle)=:(a_1, a_2, \cdots, a_{m_1+1}), \\
\textbf{b} &=:& (\langle A_1Y, Y\rangle, \langle A_2Y, Y\rangle,\cdots, \langle A_{m_1+1}Y, Y\rangle).
\end{eqnarray*}
Notice that $\widetilde{A}=\langle \textbf{a}, \textbf{b}\rangle\leq
|\textbf{a}|\cdot|\textbf{b}|\leq
\frac{1}{2}(|\textbf{a}|^2+|\textbf{b}|^2)$, and
\begin{eqnarray*}
\widetilde{A}=\frac{1}{2}(|\textbf{a}|^2+|\textbf{b}|^2) 
  &\Longleftrightarrow& \textbf{a}=\textbf{b}
\end{eqnarray*}

Define a function on the unit tangent bundle by
\begin{eqnarray*}
\phi: S(TM_1) &\longrightarrow& \mathbb{R} \\
  X &\longmapsto& \sum_{\alpha=1}^{m_1+1}\langle A_{\alpha}X, X\rangle^2
\end{eqnarray*}
For any curve $X(t)$ in $S(TM_1)$ with $X(0)=X$ a maximum point, we
have
\begin{equation*}
0 = \frac{d}{dt}\big|_{t=0}~\phi(X(t))
  = 4~\langle \sum_{\alpha=1}^{m_1+1}\langle A_{\alpha}X, X\rangle A_{\alpha}X, X'(0)~\rangle, 
\end{equation*}
which implies that
$$\sum_{\alpha=1}^{m_1+1}\langle A_{\alpha}X, X\rangle A_{\alpha}X=c\cdot X,$$
for some number $c$. Hence
$$\sum_{\alpha=1}^{m_1+1}\langle A_{\alpha}X, X\rangle^2=c \langle X, X\rangle=c\geq 0.$$

If $c=0$, then $\widetilde{A}=0$. Thus we are left to consider
$c>0$. For any orthonormal normal vectors $\{N_1, N_2, \cdots,
N_{m_1+1}\}$ of $M_1$ in the unit sphere, we denote a unit normal
vector by $\displaystyle
N=:\frac{1}{|\textbf{a}|}\sum_{\alpha=1}^{m_1+1}a_{\alpha}N_{\alpha}=\frac{1}{\sqrt{c}}\sum_{\alpha=1}^{m_1+1}a_{\alpha}N_{\alpha}$.
Then it is clear that $A_{N}X=\sqrt{c}\cdot X$.

On the other hand, recall that for any unit tangent vector on a
focal submanifold with $g=4$, the corresponding principal curvatures
are $\pm 1$ and $0$ (cf. \cite{CR}). Thus $c=1$, which leads
directly that $\widetilde{A}\leq 1$. \hfill $\Box$ \vspace{2mm}

Now we continue proving Proposition 1.1.
Combining Lemma \ref{lem Rie prod} with the fact $\widetilde{B}\geq
0$, we can conclude that
\begin{equation}\label{5.1}
Sec\leq 2.
\end{equation}


Recall that for an orthonormal basis $\{X=:e_1, e_2,\cdots, e_{m_1+2m_2}\}$ of $T_xM_1$,
the Gauss equation leads the Ricci curvature $\rho(X)$ to be
\begin{eqnarray}\label{Ricci Rie Prod}
\rho(X) 
   &=& \sum_{i=2}^{m_1+2m_2}\{1+\sum_{\alpha=1}^{m_1+1}\langle A_{\alpha}X, X\rangle\langle A_{\alpha}e_i, e_i\rangle-\sum_{\alpha=1}^{m_1+1}\langle A_{\alpha}X, e_i\rangle^2\} \\
   &=& m_1+2m_2-1-\sum_{\alpha=1}^{m_1+1}|A_{\alpha}X|^2\nonumber\\
   &\geq& 2(m_2-1),\nonumber
\end{eqnarray}
since $A_{\alpha}$ $(\alpha=1,\cdots, m_1+1)$ is trace free with
eigenvalues $\pm 1$ and $0$.

Suppose $M_1$ is a Riemannian product. Using the K\"{u}nneth formula and the Betti
numbers of the focal submanifolds given in the proof of Lemma 3.1, we can assert that $M_1\cong N_1^{m_2}\times N_2^{m_1+m_2}$. Thus for $X=e_1,
\cdots, e_{m_2}\in T_xN_1$, (\ref{5.1}) and (\ref{Ricci Rie Prod})
lead to
$$2(m_2-1)\leq \rho(X)=Sec(X\wedge e_2)+\cdots +Sec(X\wedge e_{m_2})\leq 2(m_2-1) $$
which implies that $\rho(X)=2(m_2-1)$, and further $Sec(X\wedge
e_i)=2$ $(i=2,\cdots, m_2)$, $Sec(e_j, e_i)=0$ $(i=1,\cdots, m_2,
j=m_2+1,\cdots, m_1+m_2)$.


However, let $\{\widetilde{\eta}_1=:\eta, \widetilde{\eta}_2, \cdots,
\widetilde{\eta}_{m_1+1}\}$ be orthonormal normal vectors at $x\in M_1$
in the unit sphere $S^{2m_1+2m_2+1}$, with respect to which, we have
$A_{\widetilde{\eta}_1}X=X$, and thus $A_{\widetilde{\eta}_{\alpha}}X=0$
for $\alpha=2, \cdots, m_1+1$. We choose a unit $Y\in T_xN_2$ with
$A_{\eta}Y\neq -Y$.
Thus $\langle Y, A_{\widetilde{\eta}_{\alpha}}X\rangle=0$
$(\alpha=1,\cdots, m_1+1)$. On these conditions, it is easy to see
that $\widetilde{A}(X, Y)\neq -1$ and $\widetilde{B}(X, Y)=0$, thus
$Sec(X\wedge Y)\neq 0$, a contradiction to the product splitting.

The proof of Proposition \ref{Rie product} is now complete. \hfill
$\Box$

\vspace{3mm}

To illustrate our examples to the open problem of Besse, it suffices
to prove Proposition \ref{locally homogeneous}. \vspace{2mm}

\noindent
\textbf{Proposition 1.2.} The focal submanifolds $M_1$ of OT-FKM type with $(m_1, m_2)=(3, 4k)$
are not intrinsically homogeneous.

\vspace{2mm}

\noindent \emph{Proof}. According to Theorem 5.1 in \cite{FKM}, the space
$\Omega$ defined by
\begin{equation*}
\{x\in M_1~|~\mathrm{there ~exists ~an ~orthonormal} ~Q_0,\cdots,
Q_3\in \Sigma(P_0,\cdots,P_3)~\mathrm{with}~Q_0\cdots Q_3x=x\}
\end{equation*}
can be expressed as
$$\{x\in M_1~|~\mathrm{there~exists~orthonormal} ~N_0,\cdots, N_3\in T_x^{\perp}M_1~\mathrm{with} ~\dim (\bigcap_{i=0}^3 Ker A_{N_i})\geq 3\}.$$
By Theorem 5.2 in \cite{FKM}, when $m_1=3$, $\Omega$ is non-empty
and $\Omega\neq M_1$. Comparing with Theorem 5.8 in \cite{FKM},
which states that $ \dim(\bigcap_{i=0}^3 Ker A_{N_i})\leq 3$ when
$m_1=3$, we can conclude that $\dim (\bigcap_{i=0}^3 Ker
A_{N_i})=3$. Thus for any $x\in \Omega$, and a unit $X\in T_xM_1$,
the Ricci curvature $\rho(X)$ takes the maximum $2l-6$ at the
$3$-dimensional subspace $\bigcap_{i=0}^3 Ker A_{N_i}$ of $T_xM_1$.

On the other hand, at any $y\in M_1\backslash \Omega$, the Ricci
curvature is less than the maximum $2l-6$.

The proof of Proposition \ref{locally homogeneous} is now complete.

\begin{ack}
The authors would like to thank the referees for very helpful comments.
The project was partially supported by the NSFC ( No. 11331002, and No. 11301027 ), the SRFDP (No. 20130003120008), the BJNSF (No. 1144013), the FRFCU (No. 2012CXQT09) and the Program for Changjiang Scholars and Innovative Research Team in University.
\end{ack}



\end{document}